\documentclass[12pt, a4paper, twoside]{article}

\usepackage[english]{babel}
\usepackage{amsmath}
\usepackage{latexsym}
\usepackage{amssymb,amsthm,amsfonts}
\usepackage[mathcal]{eucal}
\usepackage{amscd}
\usepackage[latin1]{inputenc}
\usepackage{color}
\def\cred{}
\def\cblu{}

\theoremstyle{plain}
\newtheorem{thm}{Theorem}[section]
\newtheorem{corol}[thm]{Corollary}
\newtheorem{lem}[thm]{Lemma}
\newtheorem{propos}[thm]{Proposition}
\newtheorem{esercize}{Exercize}[section]

\theoremstyle{definition}

\theoremstyle{remark}
\newtheorem{remark}[thm]{Remark}

\newenvironment{theor}{\vspace{-6mm}\begin{itemize}\item[]\begin{thm}}{\end{thm}\end{itemize}}

\newcommand{\N}{\mathbb{N}}                      
\newcommand{\R}{\mathbb{R}}                      
\newcommand{\LL}{L^2(\Omega)}                    
\newcommand{\HH}{H^1(\Omega)}                    
\newcommand{\vett}[2]{#1(0,T; #2)}               

\DeclareMathOperator{\Div}{div}                                                   
\newcommand{\abs}[1]{\left| #1 \right|}                                           
\newcommand{\norm}[1]{\left\| #1 \right\|}                                        
\newcommand{\scal}[2]{\left( #1, \, #2 \right)}                                   
\newcommand{\dual}[2]{\left\langle #1, \, #2 \right\rangle}                       

\newcommand{\nh}[1]{\norm{#1}^2_H}                                                
\newcommand{\Pab}{\left(\textbf{P}_{\alpha, \beta}\right)}
\newcommand{\Pa}{\left(\textbf{P}_{\alpha}\right)}

\numberwithin{equation}{section}

\begin{document}

\pagestyle{myheadings}
\newcommand\testopari{\sc Giacomo Canevari and Pierluigi Colli}
\newcommand\testodispari{\sc Solvability and asymptotic analysis
of a phase field system}
\markboth{\testodispari}{\testopari}
\renewcommand{\thefootnote}{\fnsymbol{footnote}}

\thispagestyle{empty}
{\cred
\begin{center}
{{\bf\huge Solvability and asymptotic analysis\\[1.5mm]
of a generalization of the Caginalp\\[3mm]
phase field system\footnote{{\bf Acknowledgment.}\quad\rm
The financial support of the MIUR-PRIN Grant 2008ZKHAHN 
\emph{``Phase transitions, hysteresis 
and multiscaling''} and of the IMATI of CNR in
Pavia is gratefully acknowledged.}}}

\vspace{6mm}
{\large\sc Giacomo Canevari and Pierluigi Colli}

\vspace{3mm}
{\sl Dipartimento di Matematica  ``F. Casorati'', Universit\`a di Pavia\\
Via Ferrata, 1, 27100 Pavia,  Italy\\
E-mail: {\tt giacomo.canevari@gmail.com \ \  pierluigi.colli@unipv.it}}

\end{center}

\vskip6mm
\begin{abstract}

We study a diffusion model of phase field type, which consists of 
a system of two partial differential equations involving as variables the thermal displacement, that is basically the time integration of temperature, and the order parameter. 
Our analysis covers the case of a non-smooth (maximal monotone) graph along with a smooth anti-monotone function in the phase equation. Thus, the 
system turns out a generalization of the well-known Caginalp phase field model for 
phase transitions when including a diffusive term for the thermal displacement in the 
balance equation. Systems of this kind have been extensively studied by Miranville 
and Quintanilla. We prove existence and uniqueness of a 
weak solution to the initial-boundary value problem, as well as
various regularity results ensuring that the solution is strong and with bounded components. Then we 
investigate the asymptotic behaviour of the solutions as the coefficient of the diffusive term for the thermal 
displacement tends to $0$ and prove convergence to the Caginalp 
phase field system as well as error 
estimates for the 
difference of the solutions.  
 
\end{abstract}

\section{Introduction} \label{intro} 
This paper is concerned with the initial and boundary value problem:
\begin{equation}
w_{tt} - \alpha\Delta w_t - \beta\Delta w + u_t = f \quad \textrm{ in } \Omega
\times (0,T)
\label{1}
\end{equation}
\begin{equation}
u_t - \Delta u + \gamma (u)  + g(u) \ni w_t \quad \textrm{ in } \Omega
\times (0,T)
\label{2}
\end{equation}
\begin{equation}
\partial_n w = \partial_n u = 0 \qquad \textrm{on } \Gamma\times (0, T)
\label{3}
\end{equation}
\begin{equation}
w(\cdot, \, 0) = w_0\, , \quad w_t (\cdot, \, 0) = v_0 \, , \quad u(\cdot, \, 0) = u_0 \qquad \textrm{ in } \Omega
\label{4}
\end{equation}
where $\Omega \subset \R^3 $ is a bounded domain with smooth boundary $\Gamma$,
$T> 0$ represents some finite time, and $\partial_n $ denotes the outward 
normal derivative on $\Gamma$. Moreover, $\alpha$ and $\beta$ are two positive 
parameters, $\gamma : \R \to 2^{\R} $ is a maximal monotone graph (one can see 
\cite[in particular pp. 43--45]{Brezis} or \cite{Barbu}), 
$g : \R \to \R$ is a Lipschitz-continuous 
function, $ f$ is a given source term in equation \eqref{1}
and $w_0,\,  v_0 , \, u_0$ stand for initial data.  The inclusion (in place of the equality) in \eqref{2}
is due to the presence of the possibly multivalued graph $\gamma$.

Equations \eqref{1}--\eqref{2} yield a system of phase field type. Such systems 
have been introduced (cf.~\cite{caginalp}) in order to include phase dissipation 
effects in the dynamics of moving interfaces arising in thermally induced 
phase transitions. In our case, we move from the following expression for the 
total free energy
\begin{equation}
\Psi (\theta, u) = \int_\Omega \left( - \frac12 \theta^2 - \theta u + \phi (u) + G(u) + \frac12 |\nabla u |^2 \right)
\label{5}
\end{equation}
where the variables  $\theta$ and $u$ denote the (relative) temperature and order parameter, 
respectively. Let us notice from the beginning that our $w$ represents the thermal 
displacement variable, related to $\theta$~by 
\begin{equation}
w(\cdot, \, t) = w_0 + (1* \theta) (\cdot, \, t) = w_0 + \int_0^t\!\! \theta (\cdot, \, s)\, 
ds , \quad \ t\in [0,T].
\label{6}
\end{equation}
In \eqref{5}, $ \phi : [0,+\infty] \to \R $ is  the convex and lower 
semicontinuous function such that $\phi(0) = 0 = \min \phi$ and its subdifferential 
$\partial \phi$ coincides with $\gamma $, while $G$ stands for a smooth, 
in general concave, function such that $G' = g$. A typical example for  $\phi $ 
and $G$ is the double obstacle case 
\begin{equation}
\phi(u) =  I_{[-1, +1]} (u) =  
\begin{cases} 0 & \text{if $|u|\leq 1$}
\\
+\infty &\text{if $|u| >1$}
\end{cases}
, \quad \ 
G(u) = 1- u^2
\label{7}
\end{equation}
so that the two wells of the sum  $\phi (u) + G(u)$ are located in $ -1$ and $+1$, 
and one of the two is preferred as minimum of the potential in \eqref{5} according to 
whether the temperature $\theta$ is negative or positive. Indeed, note the presence 
of the term  $- \theta u $ besides  $\phi (u) + G(u)$ in the expression of $\Psi$. 

The example given in \eqref{7} is inspired by the systematic approach of Michel Fr\'emond to 
non-smooth thermomechanics: we refer to the monography \cite{fremond} which also deals with the phase change models. In the case of \eqref{7} the subdifferential
 of the indicator function of the interval $[-1, +1]$ reads 
$$ \xi \in \partial I_{[-1, +1]} (u) \quad \hbox{ if and only if } \quad 
\xi \ \left\{
\begin{array}{ll}
\displaystyle
\leq \, 0 \   &\hbox{if } \ u=-1   
\\[0.1cm]
= \, 0 \   &\hbox{if } \ |u| < 1  
\\[0.1cm]
\geq \, 0 \  &\hbox{if } \  u = + 1  
\\[0.1cm]
\end{array}
\right. .
$$
Let us point out that, with a different terminology motivated by earlier studies 
on the Stefan problem~\cite{duvaut}, some authors (cf.~\cite{fremond}) prefer to name ``freezing index'' the variable $w$ defined by \eqref{6}, having also in mind applications to frost propagation in porous media. 

Another meaningful variable of the Stefan problem is the enthalpy $e$, which in our case is defined by
$$ 
e= - d_\theta \Psi \quad (- \hbox{ the variational derivative of $\Psi $ with respect to } \theta ) ,
$$    
whence $ e = \theta + u = w_t + u $. Then, the governing balance and phase equations 
are given~by 
\begin{equation}
e_{t} + \Div {\bf q}  = f 
\label{1phys}
\end{equation}
\begin{equation}
u_t +  d_u \Psi =0
\label{2phys}
\end{equation}
where ${\bf q} $ denotes the thermal flux vector and $d_u \Psi$ stands for the 
variational derivative of $\Psi$ with  respect to $u$. Hence, \eqref{2phys} reduces 
exactly to \eqref{2} along with the Neumann homogeneous boundary condition for $u$. 
If we assume the classical Fourier law $ {\bf q} = - \nabla \theta $ (for the moment let us take the heat 
conductivity coefficient just equal to 1), 
then  \eqref{1phys} is nothing but the usual energy balance equation as in the 
Caginalp model~\cite{caginalp}. This is also as in the weak formulation 
of the Stefan problem, in which the mere pointwise inclusion $ u \in \left( \partial I_{[-1, +1]}\right)^{-1} ( \theta)$, or equivalently $ \theta \in  \partial I_{[-1, +1]} ( u)$, replaces \eqref{2}. 

Another approach, which is by now well established, consists in adopting the so-called Cattaneo-Maxwell law (see, e.g., \cite{cgg1, MQ1} and references therein): 
such a law reads 
\begin{equation}
{\bf q} + \varepsilon {\bf q}_t = - \nabla \theta , \quad \hbox{ for } \, \varepsilon 
> 0 \, \hbox{ small},
\label{qeps}
\end{equation} 
and leads to the following equation 
\[
\varepsilon \theta_{tt} + \theta - \Delta \theta + \varepsilon u_{tt} +  u_t = f \quad \textrm{ in } \Omega
\times (0,T)
\label{1mv}
\]
which has been investigated in  \cite{MQ1}. On the other hand, if we solve \eqref{qeps} with respect to ${\bf q} $ we find
$$
{\bf q} = {\bf q_0} + k* \nabla \theta , \ \hbox{ where } \, (k* \nabla \theta) 
(x,t)  := \int_0^t  \!\! k(t-s) \nabla \theta (x,s) ds ,
$$
${\bf q_0} (x,t)$ is known and can be incorporated in the source term, $k (t) $
is a given kernel (depending on $\varepsilon$ of course): from \eqref{1phys}
we obtain the balance equation for the standard phase field model 
with memory which has a hyperbolic character and has been extensively studied in 
\cite{cgg1, cgg2}.

In \cite{gn1,gn2,gn3,gn4} Green and Naghdi presented an alternative approach based on  
a thermomechanical theory of deformable media. This theory takes advantage of an  
entropy balance rather than the usual entropy inequality. If we restrict our attention to the heat conduction, these authors
proposed three different theories, labeled as type I, type II and type III, 
respectively. In particular, when type I is linearized, we recover the classical 
theory based on the Fourier law
\begin{equation}
{\bf q} = - \alpha \nabla w_t  , \quad \alpha >0 \ \hbox{ (type I). }  
\label{typeI}
\end{equation} 
Furthermore, linearized versions of the two other theories yield
\begin{equation}
{\bf q} = - \beta \nabla w , \quad \beta >0 \ \hbox{ (type II) }  
\label{typeII}
\end{equation} 
and 
\begin{equation}
{\bf q} = - \alpha \nabla w_t - \beta \nabla w   \quad \hbox{(type III). }  
\label{typeIII}
\end{equation} 
Note that here we have used the thermal displacement \eqref{6} (instead of $\theta$)
to write such laws. We also point out that \eqref{typeII}--\eqref{typeIII}  
have been recently discussed, applied and compared by 
Miranville and Quintanilla in \cite{MQ2, MQ3, MQ4} (there the reader can find a rich 
list of references as well). In particular, \eqref{typeIII} leads via 
\eqref{1phys} to our equation \eqref{1}; further, a no flux boundary condition for $ {\bf 
q}$ corresponds to $ \partial_n w = 0  $ in \eqref{3}.

Thus, the system \eqref{1}--\eqref{4} results from \eqref{1phys}--\eqref{2phys} when 
\eqref{5} and \eqref{typeIII} are postulated. We are interested in the study of 
existence, uniqueness, regularity of the solution to the initial-boundary value 
problem \eqref{1}--\eqref{4} when $\gamma$ is an arbitrary maximal monotone graph, 
possibly multivalued, singular and with bounded domain. Of course,
the case of $\Psi$ shaped by a multiwell potential 
$ u \mapsto - w_t u + \phi (u) + G(u)$ is recovered as a sample. Then we study the asymptotic behaviour of the 
problem as $\beta \searrow 0$, obtaining convergence of solutions to the problem with 
$\beta=0$, which corresponds to \eqref{typeI}, the (type I) case of Green and Naghdi.
We also prove two error estimates of the difference of solutions in suitable norms,
showing a linear rate of convergence in both estimates. In a subsequent study 
we would like to address the investigation of the analogous limit $\alpha \searrow 0$
to obtain the (type II) case in \eqref{typeII}. 

The paper is organized as follows. In Section~\ref{wepo} we state the main results
related to the problem~\eqref{1}--\eqref{4}: existence and uniqueness of a weak 
solution, regularity results yielding a strong solution, further regularity results 
ensuring the boundedness of $u, \, w_t $ and of the appropriate selection of $\gamma (u)$. 
Section~\ref{sec: Pa} contains the related statements. Then we 
investigate the asymptotic limit as 
$\beta \searrow 0$: precisely, the convergence result and the error estimates under
different assumptions on the data. In Section~\ref{no-un} we introduce some notation 
and present the uniqueness proof. The approximation of the 
problem~\eqref{1}--\eqref{4} via a Faedo-Galerkin scheme and the derivation of the
uniform a priori estimates are carried out in Section~\ref{app}. Regularity and 
boundedness properties for the solutions are proved in Sections~\ref{reg1}--\ref{reg3}. 
Finally, the details of the asymptotic analysis as $\beta \searrow 0$ are developed 
in Section~\ref{beta=0}.

}%
\section{Well-posedness and regularity for $\alpha, \beta >0$}
\label{wepo}
We point out the assumptions on the data and state clearly 
the formulation of the problem and the main results we achieve.
Let $\Omega\subseteq\R^3$ be a bounded {\cred smooth domain}
with boundary $\Gamma = \partial\Omega$ {\cred and} let $T > 0$. Set {\cred $Q:= \Omega \times (0,T)$. We assume that}
\begin{equation}
\alpha\, , \; \beta \in (0, +\infty)
\label{alpha, beta}
\end{equation}
\begin{equation}
\label{f}
f \in \vett{L^2}{\HH'} + \vett{L^1}{\LL}
\end{equation}
\begin{equation}
\gamma\subseteq\R\times\R \, \textrm{ is a maximal monotone graph, with } \, \gamma(0)\ni 0
\label{gamma}
\end{equation}
\begin{equation}
\phi:\R\longrightarrow [0, +\infty] \, \textrm{ is convex and lower-semicontinuous}
\label{phi}
\end{equation}
\begin{equation}
\phi(0) = 0 \ \textrm{ and } \ \partial\phi = \gamma
\label{phi gamma}
\end{equation}
\begin{equation}
g: \R \longrightarrow \R \, \textrm{ is Lipschitz-continuous}
\label{g}
\end{equation}
\begin{equation}
w_0 \in \HH \, , \quad v_0\in\LL \, , \quad u_0\in\LL \, , \quad \phi(u_0)\in L^1(\Omega) .
\label{initial data}
\end{equation}
The effective domain of $\gamma$ will be denoted by $D(\gamma)$. We consider

\bigskip 
\noindent
\textbf{Problem $\Pab$.} Find $(w, \, u, \, \xi)$ satisfying
\begin{equation}
w \in \vett{W^{1, \, \infty}}{\LL} \cap \vett{H^1}{\HH}
\label{w}
\end{equation}
\begin{equation}
w_{tt} \in \vett{L^1}{\HH'}
\label{w_tt}
\end{equation}
\begin{equation}
u \in \vett{H^1}{\HH'} \cap C^0\left([0, T]; \, \LL\right) \cap \vett{L^2}{\HH}
\label{u}
\end{equation}
\begin{equation}
\xi\in L^2(Q) \, , \qquad u\in D(\gamma) \ \textrm{ and } \ \xi\in\gamma(u)\  \textrm{ a.e. in } \, Q
\label{xi}
\end{equation}
\begin{equation}
\begin{split}
\dual{w_{tt}(t)}{v} + \alpha\scal{\nabla w_t(t)}{\nabla v}_{\LL} + \beta\scal{\nabla w(t)}{\nabla v}_{\LL}
+ \dual{u_t(t)}{v} = \dual{f(t)}{v} \\ \qquad\textrm{for all } v\in\HH \textrm{ and a.a. } t\in (0, T)
\end{split}
\label{eq. A weak}
\end{equation}
\begin{equation}
\begin{split}
\dual{u_t(t)}{v} + \scal{\nabla u(t)}{\nabla v}_{\LL} + \scal{\xi(t)}{v}_{\LL} + \scal{g(u)(t)}{v}_{\LL}
= \scal{w_t(t)}{v}_{\LL} \\ \qquad \textrm{for all } v\in\HH \textrm{ and a.a. } t\in(0, T)
\end{split}
\label{eq. B weak}
\end{equation}
\begin{equation}
{\cred
w(0) = w_0\ \textrm{ in } \, H^1(\Omega) \, , \quad w_t(0) = v_0 \ \textrm{ in } \, H^1(\Omega)'\, , \quad u(0) = u_0 \ \textrm{ in } \, L^2(\Omega). }
\label{initial condition}
\end{equation}
We can prove the well-posedness of this problem.

\begin{theor}[Existence and uniqueness] Let assumptions \eqref{alpha, beta}--\eqref{initial data} hold. Then Problem~$\Pab$ has a unique solution.
\end{theor}

Next, in addition to \eqref{alpha, beta}--\eqref{initial data}, we suppose
\begin{equation}
f \in L^2(0,T;L^2(\Omega)) + \vett{L^1}{\HH}
\label{f strong}
\end{equation}
\begin{equation}
w_0 \in H^2(\Omega) \, , \quad \partial_n w_0 = 0 \, \textrm{ on } \, {\cred \Gamma} \, , \quad  v_0 \in \HH \, , \quad u_0\in \HH \, ;
\label{initial data strong}
\end{equation}
in this case, we are able to prove a regularity result, which allows us to {\cred solve a strong formulation
of Problem~$\Pab$}.

\begin{theor}[Regularity and strong solution]
\label{th: strong solution}
Assume \eqref{f strong}--\eqref{initial data strong} in addition to \eqref{alpha, beta}--\eqref{initial data}. Then the unique solution $(w, \, u, \, \xi)$ of Problem~$\Pab$ fulfills
\begin{equation}
w \in \vett{W^{1, \, \infty}}{\HH} \cap \vett{H^1}{H^2(\Omega)}
\label{w strong}
\end{equation}
\begin{equation}
w_{tt} \in \vett{L^1}{\LL} 
\label{w_tt strong}
\end{equation}
\begin{equation}
u \in \vett{H^1}{\LL} \cap C^0\left([0, T]; \, \HH\right) \cap \vett{L^2}{H^2(\Omega)} \, .
\label{u strong}
\end{equation}
In particular, $(w, \, u, \, \xi)$ solves Problem~$\Pab$ in a strong sense, that is,
{\cred $w$ and $u$ satisfy}
\[
w_{tt} - \alpha\Delta w_t - \beta\Delta w + u_t = f \quad \textrm{ a.e. in } Q
\]
\[
u_t - \Delta u + \xi + g(u) = w_t, \quad {\cred \xi \in \gamma(u)}  \quad \textrm{ a.e. in } Q
\]
\[
\partial_n w = \partial_n u = 0 \quad \textrm{ a.e. on } \Gamma\times (0, T) \, .
\]
\end{theor}

The aim of the subsequent results is to provide $L^\infty$ estimates. We will need to 
strengthen again the hypotheses on the initial data. {\cred For $s\in D(\gamma)$ let 
us denote} by $\gamma^0(s)$  the element of $\gamma(s)$ having minimal modulus.  
{\cred Then, we require that}
\begin{equation}
u_0 \in H^2(\Omega) \, , \quad \partial_n u_0 = 0 \,\textrm{ on } \, {\cred \Gamma}
\label{u0 strong}
\end{equation}
\begin{equation}
u_0 \in D(\gamma) \quad \textrm{a.e. in } \Omega \, , \quad \ \gamma^0(u_0)\in\LL \, .
\label{u0 gamma}
\end{equation}

\begin{theor}[Further regularity] \label{th: regularity}
If the conditions \eqref{alpha, beta}--\eqref{initial data}, \eqref{f strong}--\eqref{initial data strong} and \eqref{u0 strong}--\eqref{u0 gamma} hold, then the solution $(w, \, u, \, \xi)$ of Problem~$\Pab$ fulfills
\begin{equation}
u \in \vett{W^{1, \, \infty}}{\LL} \cap \vett{H^1}{\HH} \cap \vett{L^\infty}{H^2(\Omega)} \, .
\label{u stronger}
\end{equation}
\end{theor}

The above results still hold if the dimension $N$ of the domain $\Omega$ is arbitary. On the other 
hand, {\cred since \eqref{u stronger} implies in particular that 
$u$ is  continuous from $[0,T]$ to the space $H^s(\Omega)$   
for all $s<2$, then, if we let $N \leq 3$ and $s$ sufficiently large, it turns out that 
$H^s(\Omega) \subset C^0(\overline\Omega) $ and consequently}
\[
u \in C^0(\overline Q) \, .
\]

Finally, we assume for the data {\cred enough regularity to get 
$L^\infty$ estimates for} $w_t$ and $\xi$. 
The hypothesis $N\leq 3$ is essential in the proof of the following {\cred result}.
\begin{theor}[$L^\infty$ estimate for $w_t$ and $\xi$] \label{th: L^infty estimate}
In addition to {\cred assumptions} \eqref{alpha, beta}--\eqref{initial data}, 
\eqref{f strong}--\eqref{initial data strong} and \eqref{u0 strong}--\eqref{u0 gamma}, we ask
\begin{equation}
f \in \vett{L^\infty}{\LL} + \vett{L^r}{\HH} \quad \textrm{ for some } \, r > 4/3 
\label{f stronger}
\end{equation}
\begin{equation}
\gamma^0(u_0) \in L^\infty(\Omega) \, .
\label{gamma strong}
\end{equation}
Then we have
\[
w_t \in L^\infty(Q) \, , \qquad \xi \in L^\infty(Q) \, .
\]
\end{theor}

{\cblu \begin{remark}
All the statements contained in this paper still hold if $\Omega\subseteq\R^3$ is, for instance, a convex polyhedron, for which standard results on Sobolev embeddings and regularity for elliptic problems apply.
\end{remark}}

\section{Asymptotic behaviour as $\beta\searrow 0$} \label{sec: Pa}
Let us fix the parameter $\alpha$ once and for all. We shall concentrate 
on the asymptotic behaviour of the solution {\cred as} $\beta \searrow 0$, 
so we let $\beta$ vary in a bounded subset of $(0, +\infty)$. We {\cred allow} 
the source term and the initial data in Problem~$\Pab$ {\cred to} vary with $\beta$, 
by replacing $f$, $w_0$, $v_0$ and $u_0$ in {\cred \eqref{eq. A weak} 
and \eqref{initial condition}} with $f_\beta$, $w_{0, \beta}$, $v_{0,\,\beta}$ 
and $u_{0, \beta}$ respectively. We will denote by $(w_\beta, \, u_\beta, \, 
\xi_\beta)$ the solution to Problem~$\Pab$.

If we set $\beta = 0$ in the statement of Problem~$\Pab$, we get a first-order system of differential equations, with respect to time, in the variable $w_t$, which is of physical relevance {\cred (recall that $w_t=\theta$)}. Anyway, we avoid this change of variable, in order to preserve the formalism. We {\cred introduce} the formulation of Problem~$\Pa$, in which $\beta$ is set to be zero.

\bigskip 
\noindent
\textbf{Problem $\Pa$.} Find $(w, \, u, \, \xi)$ satisfying \eqref{w}--\eqref{xi} as well as {\cred
\begin{equation}
\begin{split}
\dual{w_{tt}(t)}{v} + \alpha\scal{\nabla w_t(t)}{\nabla v}_{\LL} + \dual{u_t(t)}{v} = \dual{f(t)}{v}\\
\hbox{for all \,$v\in\HH$ \, and a.a. \,$t\in (0,T)$} 
\end{split}
\label{eq. Aa weak}
\end{equation}
\begin{equation}
\begin{split}
\dual{u_t(t)}{v} + \scal{\nabla u(t)}{\nabla v}_{\LL} + \scal{(\xi + g(u))(t)}{v}_{\LL} = \scal{w_t(t)}{v}_{\LL}\\
\hbox{for all \, $v\in\HH$ \, and a.a.\, $t\in (0,T)$}
\end{split}
\label{eq. Ba weak}
\end{equation}
\begin{equation}
{\cred
w(0) = w_0\ \textrm{ in } \, H^1(\Omega) \, , \quad w_t(0) = v_0 \ \textrm{ in } \, H^1(\Omega)'\, , \quad u(0) = u_0 \ \textrm{ in } \, L^2(\Omega). }
\label{initial condition a}
\end{equation}
}%

We state at first the well-posedness of Problem~$\Pa$ and a convergence result.

\begin{theor}[Well-posedness for $\Pa$] \label{th: well-posedness Pa} If the hypotheses \eqref{f}--\eqref{initial data} hold, then Problem~$\Pa$ admits exactly one solution.
\end{theor}
\begin{theor}[Convergence as $\beta\searrow 0$] \label{th: convergence}
We assume {\cred \eqref{f}--\eqref{initial data} and}
\begin{equation}
f_\beta \rightharpoonup f \quad \textrm{in } \vett{L^2}{\HH'} + \vett{L^1}{\LL}
\label{f convergence beta}
\end{equation}
\begin{equation}
w_{0, \beta} \rightharpoonup w_0 \quad \textrm{in } \HH \, , \qquad v_{0, \beta} \rightharpoonup v_0 \, , \quad u_{0, \beta} \rightharpoonup u_0 \quad \textrm{in } \LL .
\label{data convergence beta}
\end{equation}
Then, the convergences
\[
w_\beta \rightharpoonup^* w \quad \textrm{in } \vett{W^{1, \, \infty}}{\LL} \, , \qquad w_\beta \rightharpoonup w \quad \textrm{in } \vett{H^1}{\HH}
\]
\[
u_\beta \rightharpoonup u \quad \textrm{in } \vett{H^1}{\HH'} \cap \vett{L^2}{\HH}
\]
\[
\xi_\beta \rightharpoonup \xi \quad \textrm{in } L^2(Q) \, .
\]
hold{\cred , where $(w,u,\xi)$ denotes the solution to  Problem  $\Pa$.}
\end{theor}

With slightly strengthened hypotheses, we are able to prove 
the strong convergence for the {\cred solution and even} to give an estimate for the convergence 
{\cred rate}.

\begin{theor}[{\cred First error estimate}] \label{th: first estimate error}
In addition to \eqref{gamma}--\eqref{g} and \eqref{f convergence beta}--\eqref{data convergence beta}, we assume
\begin{equation}
\norm{f_\beta - f}_{\vett{L^2}{\HH'} + \vett{L^1}{\LL}} \leq c\,\beta
\label{f rate a}
\end{equation}
\begin{equation}
{\cred \norm{w_{0, \beta} - w_0}_{H^1(\Omega)} }+ \norm{v_{0, \beta} - v_0}_{H^1(\Omega)'} + \norm{u_{0, \beta} - u_0}_{\LL} \leq c\, \beta 
\label{data rate a}
\end{equation}
for some constant $c$ which is independent of $\beta$. Then {\cred one has} the estimate{\cred
\begin{equation}
\begin{split}
\norm{w_\beta - w}_{\vett{H^1}{\LL}\cap\vett{L^\infty}{H^1(\Omega)}} 
\hskip2cm\\
+ \norm{u_\beta - u}_{\vett{L^\infty}{\LL}\cap\vett{L^2}{H^1(\Omega)}} \leq c\,\beta
\end{split}
\label{stimaerr1}
\end{equation}
}%
where $c$ does not depend on $\beta$.
\end{theor}

If $\gamma$ is a (single-valued) smooth function, and if enough regularity on the data is assumed, it is possible to obtain much stronger estimates. The assumption $N\leq 3$ on the spatial dimension is essential for the proof of the following result.

\begin{theor}[{\cred Second error estimate}] 
\label{th: second estimate error}
{\cred Let \eqref{gamma}--\eqref{g}, \eqref{f convergence beta}--\eqref{data convergence beta} hold {\cred and} 
\begin{equation}
\gamma: D(\gamma)\longrightarrow \R \ \textrm{ be single-valued and locally Lipschitz-continuous.}
\label{gamma lipshitz}
\end{equation}
Moreover, assume that the data 
$\{ f_\beta, \, w_{0, \beta},\, v_{0, \beta}, \, u_{0, \beta} \}$, 
as well as $\{ f, \, w_{0},\, v_{0}, \, u_{0} \}$, satisfy \eqref{f strong}--\eqref{initial data strong}, \eqref{u0 strong}--\eqref{u0 gamma}, 
\eqref{f stronger}--\eqref{gamma strong} along with 
\begin{equation}
\norm{f_\beta}_{\vett{L^\infty}{L^2(\Omega)} + \vett{L^r}{H^1(\Omega)} }
+ \norm{ u_{0, \beta}}_{H^2(\Omega)} 
+ \norm{ \gamma(u_{0, \beta}) }_{L^\infty (\Omega)} \leq c
\label{debole1}
\end{equation}
\begin{equation}
\norm{f_\beta - f}_{\vett{L^2}{L^2(\Omega)} + \vett{L^1}{H^1(\Omega)}} \leq c\, \beta
\label{f rate b}
\end{equation}
\begin{equation}
\norm{w_{0, \beta} - w_0}_{H^2(\Omega)} + \norm{v_{0, \beta} - v_0}_{\HH} + \norm{u_{0, \beta} - u_0}_{\HH} \leq c\, \beta 
\label{data rate strong}
\end{equation}
where $r>4/3$. Then the estimate
\begin{equation}
\begin{split}
\norm{w_\beta - w}_{\vett{W^{1,\infty}}{H^1(\Omega)}\cap\vett{H^1}{H^2(\Omega)}} 
\hskip4cm\\ 
+ \norm{u_\beta - u}_{\vett{H^1}{L^2(\Omega)}\cap\vett{L^\infty}{H^1(\Omega)}\cap\vett{L^2}{H^2(\Omega)}} \leq c \, \beta
\end{split} 
\label{stimaerr2}
\end{equation}
holds for a suitable constant $c$, which may depend on $\alpha$ but not on $\beta$.
}%
\end{theor}

{\cred
\section{Notation and uniqueness proof}
\label{no-un}
Before facing the proof of all the 
}%
results, for the sake of convenience we fix some notation:
\[ Q_t = \Omega \times {\cred (0, t)  \quad \textrm{for } 0 \leq t \leq T , \quad \ Q=Q_T ,}  \]
\[ H = \LL \, , \ \quad V = \HH \, , \quad\ 
{\cred W = \left\{v\in H^2(\Omega): \: \partial_n v = 0 \quad \textrm{a.e. on } \Gamma \right\}. }
\]
We embed $H$ in $V'$, by means of the formula
\[ \dual{y}{v} = \scal{y}{v}_H \qquad \textrm{for all } y\in H \, , \; v\in\ V \, . \]
Furthermore, the same symbol $\norm{\cdot}_H$ will denote both the norm 
in $\LL$ and in $\LL^{{\cred N}}$; we behave similarly with $\norm{\cdot}_V$. 
If $a$, $b$ are functions of space and time variables, we introduce the 
convolution product with respect to time
\[ 
(a * b)(t) = \int_0^t a(s)b(t - s) ds \, , \qquad 0 \leq t \leq T \, .
\]
We also point out that the symbols $c$, $c_i$ -- even in the same formula -- stand for 
different constants, depending on $\Omega$, $T$ and the data, but not on the parameters 
$\alpha$, $\beta$. However, as we will be interested in the study of convergence as 
$\beta\searrow 0$,  if a constant $c$ depends on $\alpha$, $\beta$ in such a way that $c$ is 
bounded whenever $\alpha$, $\beta$ lie bounded, then we will accept the notation $c$. A constant depending on the data and 
on $\alpha$, but not on $\beta$, may be denoted by $c_\alpha$ or $c_{\alpha, i}$ {\cred or simply $c$, 
as it will happen in Section~\ref{beta=0}.} 

{\cred In our computations, we will often exploit the H\"older and Young inequalities 
to infer}
\[
\int_{Q_t} ab \leq \frac{1}{2\sigma} \int_0^t \nh{a(s)}ds + \frac{\sigma}{2} \int_0^t 
\nh{b(s)} ds 
\]
where $a, \, b \in L^2(Q)$ and $\sigma > 0$ is arbitrary. We point out another inequality 
which will turn out to be useful: if $\varphi\in\vett{H^1}{H}$, then the fundamental {\cred 
t}heorem of calculus and the H\"older inequality entail
\begin{equation}
{\cred \nh{\varphi(t)} =  \left\| \varphi(0) + \int_0^t\varphi_t(s)ds\right\|_H^2 
\leq 2\nh{\varphi(0)} + 2T \int_0^t \nh{\varphi_t(s)} ds  } 
\label{fond. t. calculus}
\end{equation}
for all $0 \leq t \leq T$. Now, let us concentrate on the uniqueness proof. 

Let $(w_1, \, u_1, \, \xi_1)$ and $(w_2, \, u_2, \, \xi_2)$ be solutions to the Problem~$\Pab$; we claim that they coincide. Setting $w = w_1 - w_2$, \mbox{$u = u_1 - u_2$} and $\xi = \xi_1 - \xi_2$, we easily get
\begin{equation}
\dual{w_{tt}(t)}{v} + \alpha\scal{\nabla w_t(t)}{\nabla v}_H + \beta\scal{\nabla w(t)}{\nabla v}_H + \dual{u_t(t)}{v} = 0
\label{A uniq.}
\end{equation}
\begin{equation}
\begin{split}
\dual{u_t(t)}{v} + \scal{\nabla u(t)}{\nabla v}_H + \scal{\xi(t)}{v}_H + \scal{g(u_1)(t) - g(u_2)(t)}{v}_H  = \scal{w_t(t)}{v}_H 
\end{split}
\label{B uniq.}
\end{equation}
for all $v\in V$ and a.a. $0 \leq t \leq T$, along with the initial conditions {\cred 
\begin{equation}
w(0) = w_t(0) = u(0) = 0 \, . \label{3-2bis}
\end{equation}
}%
We choose $v = u(t)$ in equation \eqref{B uniq.} and integrate over $(0, t)$; thus, we obtain
\[
\frac{1}{2}\nh{u(t)} + \int_0^t \nh{\nabla u(s)} ds + \int_{Q_t}\xi u = - \int_{Q_t} \left(g(u_1) - g(u_2)\right)u + \int_{Q_t} w_t u \, .
\]
Accounting for the Lipschitz-continuity of $g$, the H\"older inequality and the monotonicity of $\gamma$, frow the above equality we easily derive
\begin{equation}
\frac{1}{2}\nh{u(t)} + \int_0^t \nh{\nabla u(s)} ds \leq c \int_0^t\nh{u(s)}ds + \int_{Q_t} w_t u \, .
\label{B uniq..}
\end{equation}
Integrating in time the equation \eqref{A uniq.} (this is possible thanks to \eqref{w}) and taking the initial data {\cred \eqref{3-2bis}} into account, we have
\begin{equation}
\scal{w_{t}(t)}{v}_H + \alpha\scal{\nabla w(t)}{\nabla v}_H + \beta\scal{1*\nabla w(t)}{\nabla v}_H + \scal{u(t)}{v}_H = 0 \, ;
\label{A new uniq.}
\end{equation}
we choose $v = w_t(t)$ in \eqref{A new uniq.} and integrate over $(0, t)$. Noticing that the equality
\begin{equation}
\scal{1*\nabla w(t)}{ {\cred \nabla} w_t(t)}_H = \frac{d}{dt} \scal{1*\nabla w(t)}{\nabla w(t)}_H - \nh{\nabla w(t)}
\label{derivative}
\end{equation}
holds, we get
\begin{equation}
\begin{split}
\int_0^t\nh{w_t(s)}ds + \frac{\alpha}{2} \nh{\nabla w(t)} = - \beta\scal{1*\nabla w(t)}{\nabla w(t)}_H \\
+ \beta\int_0^t\nh{\nabla w(s)} ds - \int_{Q_t} uw_t \, .
\end{split}
\label{A new uniq..}
\end{equation}
The H\"older inequality and \eqref{fond. t. calculus} allow us to deal with the right-hand side of this formula:
\begin{equation}
- \beta\scal{1*\nabla w(t)}{\nabla w(t)}_H \leq \frac{c\beta^2}{\alpha} \int_0^t\nh{\nabla w(s)}ds + \frac{\alpha}{4}\nh{\nabla w(t)} \, . \\
\label{A diseg}
\end{equation}
Collecting {\cred now} \eqref{B uniq..}, \eqref{A new uniq..} 
and \eqref{A diseg}, it follows that
\[ 
\begin{split}
\frac{1}{2}\nh{u(t)} + \int_0^t \nh{\nabla u(s)} ds + \int_0^t\nh{w_t(s)}ds + \frac{\alpha}{4} \nh{\nabla w(t)} \\
\leq c \int_0^t\nh{u(s)}ds + c\left(\beta + \frac{\beta^2}{\alpha}\right) \int_0^t\nh{\nabla w(s)}ds \, ;
\end{split}
\]
{\cred then, by} applying the Gronwall lemma {\cred and recalling \eqref{3-2bis}}, we obtain $
{\cred  u=w= 0}$ almost everywhere in $Q$. A comparison in \eqref{eq. B weak} and the density of $H^1(Q)$ as a subspace of $L^2(Q)$ entail $\xi = 0$ almost everywhere in $Q$; thus, the proof of uniqueness is complete.

\section{Approximation and a priori estimates}
\label{app}
We are going to prove the existence of a solution to Problem~$\Pab$ via a 
Faedo-Galerkin method. {\cred First, we approximate the graph $\gamma$ with its Yosida regularization: for all $\varepsilon \in (0,1]$ say, we let
\[
{\cred
\gamma_\varepsilon := \frac{1}{\varepsilon}\left\{I - \left(I + \varepsilon\gamma\right)^{-1}\right\} \quad \hbox{ and } \quad
\phi_\varepsilon(s) := \min_{\tau\in\R}\left\{\frac{1}{2\varepsilon}\abs{\tau - s}^2 + \phi(\tau)\right\} \quad  \hbox{ for }\,  s\in\R 
}
\]
where $I$ denotes the identity on $\R$. We recall that $\phi_\varepsilon$ is a nonnegative, convex and differentiable function, $\gamma_\varepsilon$ is Lipschitz-continuous, monotone and {\cred
\begin{equation}
\label{prope}
\gamma_\varepsilon(0) = 0 \, , \quad \phi_\varepsilon ' = \gamma_\varepsilon \, , \quad 0 \leq \phi_\varepsilon(s) \leq \phi(s) \, , \quad \abs{\gamma_\varepsilon(s)} \leq \abs{\gamma^0(s)} \ \quad \forall\  \varepsilon > 0, \  s\in\R
\end{equation}
}%
{\cred (see, e.g., \cite[Prop.~2.6, p.~28 and Prop.~2.11, p.39]{Brezis} or \cite[pp.~57--58]{Barbu}).}

We look for a solution of the approximating problem} 
in a finite-dimensional subspace $V_n\subseteq V$, 
chosing a sequence $\left\{V_n\right\}$ filling up $V$; then we get a priori estimates 
and use compactness arguments to take the limit {\cred as} $n \longrightarrow 
+\infty$. {\cred In a second step we let $\varepsilon \searrow 0.$}

A special choice of the approximating subspaces will be useful. Let $\left\{v_i\right\}_{i\in\N}$ be an orthonormal basis for $V$ satisfing
\begin{equation}
{\cred - \Delta v_i = \lambda_i v_i \quad \textrm{ in }\,  \Omega ,  \quad \quad
\partial_n v_i = 0  \quad \textrm{ on } \, \Gamma }
\label{numeroform}
\end{equation}
where $\left\{\lambda_i\right\}_{i\in\N}$ are the eigenvalues of the Laplace operator; also, let $V_n$ be the subspace of $V$ spanned by $v_1, \, \ldots, \, v_n$, for all $n\in\N$. Thus, we have defined an increasing sequence of subspaces, whose {\cred union} is dense in $V$, and hence in $H$; furthermore, we notice that the regularity of $\cred \Omega$ implies ${\cred V_n\subseteq W }$, for all $n\in\N$.

As approximations of the data $w_0$, $v_0$, $u_0$ we choose the projections on $V_n$: let $w_{0, n}$ be the projection of $w_0$, with respect to $V$, and let $v_{0, n}$, $u_{0, n}$ be the projections of $v_0$, $u_0$, with respect to $H$. We notice that
\begin{equation}
w_{0, n} \longrightarrow w_0 \ \textrm{ in } V \, , \quad v_{0, n} \longrightarrow v_0 \ 
\textrm{ in } H \, , \quad u_{0, n} \longrightarrow u_0 \ \textrm{ in } H \, .
\label{data convergence}
\end{equation}
We also need to regularize the source term $f$: so, we first write {\cred
\begin{equation}
f = f^{(1)} + f^{(2)} \, , \quad \textrm{where } \, f^{(1)}\in\vett{L^2}{V'} \,  \textrm{ and } \, f^{(2)}\in \vett{L^1}{H} \, ,
\label{f splittata}
\end{equation}
}%
then we assume $f_n^{(1)}$, $f_n^{(2)}$ to be functions in $C^0\left([0, T]; \, V'\right)$, $C^0\left([0, T]; \, H\right)$ respectively, such that
\begin{equation}
f_n^{(1)} \longrightarrow f^{(1)} \ \textrm{ in } \vett{L^2}{V'} \, , \quad 
f_n^{(2)} \longrightarrow f^{(2)} \ \textrm{ in } \vett{L^1}{H} \, ;
\label{f convergence}
\end{equation}
we also set $f_n = f_n^{(1)} + f_n^{(2)}$. 


Now we are ready to state the approximated problem. For the sake of simplicity, we do not specify explicitly the dependency {\cred on} $\varepsilon$ in the solution.

\textbf{Problem $\Pab_{n, \, \varepsilon}$}. Find $T_n\in (0, T]$ and $(w_n,  u_n)$ satisfying
\[
w_n\in C^2([0, T_n]; \, V_n) \, , \qquad u_n\in C^1([0, T_n]; \, V_n)
\]
\begin{equation}
\begin{split}
\scal{\partial^2_{t} w_n(t)}{v}_H + \alpha \scal{\nabla \partial_t w_n (t)}{\nabla v}_H + \beta \scal{\nabla w_n(t)}{\nabla v}_H + \scal{\partial_t u_n(t)}{v}_H\\
= \dual{f_n(t)}{v} \qquad \textrm{for all } v\in V_n \textrm{ and all } t \in [0,  T_n]
\end{split}
\label{eq. A n}
\end{equation}
\begin{equation}
\begin{split}
\scal{\partial_{t} u_n(t)}{v}_H + \scal{\nabla u_n (t)}{\nabla v}_H + \scal{\gamma_\varepsilon(u_n)(t)}{v}_H + \scal{g(u_n)(t)}{v}_H \\
= \scal{\partial_t w_n(t)}{v}_H \qquad \textrm{for all } v\in V_n \textrm{ and all } t \in [0,  T_n]
\end{split}
\label{eq. B n}
\end{equation}
{\cred
\begin{equation}
\label{in.cond}
w_n(0) = w_{0, n} \, , \qquad \partial_t w_n(0) = v_{0, n} \, , \qquad u_{n}(0) = u_{0, n} \, . 
\end{equation}
}%
Writing $w_n$ and $u_n$ as linear combinations of $v_1, \, \ldots , \, v_n$ with 
time-dependent coefficients, and testing equations \eqref{eq. A n} and \eqref{eq. B n} {\cred by} $v = v_1, \, \ldots , \, v_n$, we obtain a system of ordinary {\cred differential}
equations, for whose local existence and uniqueness standard results apply. Thus, 
Problem~$\Pab_{n, \, \varepsilon}$ admits a solution, defined on some interval $[0, 
\, T_n]$. The following estimates imply that these solutions can be extended over the 
whole interval $[0, T]$.

\bigskip
\textbf{First a priori estimate.} We choose $v = {\cred u_n}(t)$ in equation \eqref{eq. B n} and integrate over~$(0, t)$:
\[
\begin{split}
\frac{1}{2}\nh{u_n(t)} + \int_0^t \nh{\nabla u_n(s)} ds + \int_{Q_t} \gamma_\varepsilon(u_n)u_n \\
= - \int_{Q_t} g(u_n)u_n + \int_{Q_t} u_n\partial_t w_n + \frac{1}{2} \nh{u_{0, n}} \, .
\end{split}
\]
The last term in the left-hand side is non negative, because $\gamma_\varepsilon$ is {\cred increasing} and $\gamma_\varepsilon(0) = 0$; it will be ignored in the following estimates. Meanwhile, the right-hand side can be easily estimated using the Lipschitz-continuity of $g$ and \eqref{data convergence}; so we get
\begin{equation}
\frac{1}{2}\nh{u_n(t)} + \int_0^t \nh{\nabla u_n(s)} ds \leq c \int_0^t \nh{u_n(s)} ds + \int_{Q_t} u_n\partial_t w_n + c \, .
\label{test 1, 2}
\end{equation}
Following the same computation {\cred as in} the uniqueness proof, we integrate equation \eqref{eq. A n} with respect to time:
\begin{equation}
\begin{split}
\scal{\partial_{t} w_n(t)}{v}_H + \alpha \scal{\nabla w_n (t)}{\nabla v}_H + \beta \scal{1*\nabla w_n(t)}{\nabla v}_H + \scal{ u_n(t)}{v}_H \qquad \\
= \dual{1*f_n^{(1)}(t)}{v} +
\scal{1*f_n^{(2)}(t)}{v}_H  + \scal{v_{0, n} + u_{0, n}}{v}_H + \alpha\scal{\nabla w_{0, n}}{\nabla v}_H 
\end{split}
\label{new A n}
\end{equation}
for all $v\in V_n$ and $0 \leq t \leq T_n$. We take $v = \partial_t w_n(t)$ in the previous equation and integrate over $(0, t)$. Recalling the identity \eqref{derivative}, we have
\begin{equation}
\begin{split}
\int_0^t \nh{\partial_tw_n(s)} ds + \frac{\alpha}{2}\nh{\nabla w_n(t)} = \sum_{i = 1}^7 T_i(t) + \frac{\alpha}{2}\nh{\nabla w_{0, n}} 
\end{split}
\label{test 1, 1}
\end{equation}
where we have set
\[
T_1(t) = \beta\int_0^t\nh{\nabla w_n(s)}ds \, , \quad T_2(t) = - \beta\scal{1*\nabla w_n(t)}{\nabla w_n(t)}_H 
\]
\[
T_3(t) = - \int_{Q_t}\!\! u_n\partial_t w_n \, , \quad 
{\cred T_4(t) = \int_0^t\!\! \left\langle 1*f_n^{(1)}(s), 
\partial_t w_n (s)\right\rangle  ds }\, , \quad T_5(t) = \int_{Q_t}\!\!\left(1*f_n^{(2)}\right)\partial_t w_n 
\]
\[
{\cred T_6(t) = \int_0^t\scal{v_{0, n} + u_{0, n}}{\partial_t w_n(s)}_H ds \, , \quad 
T_7(t) = \alpha \int_0^t\scal{\nabla w_{0, n}}{\nabla\partial_t w_n(s)}_H ds }
\, .
\]
We do not need any estimate on terms $T_1$ and $T_3$. With simple applications of the H\"older inequality, we estimate $T_2$, $T_5$ and $T_6$:
\[
T_2(t) \leq \frac{\alpha}{8}\nh{\nabla w_n(t)} + \frac{c\beta^2}{\alpha}\int_0^t\nh{\nabla w_n(s)} ds
\]
\[
T_5(t) \leq \frac{1}{4}\int_0^t \nh{\partial_t w_n(s)} ds + \int_0^t \nh{1*f^{(2)}_n(s)} ds
\]
\[
T_6(t) \leq \frac{1}{4}\int_0^t \nh{\partial_tw_n(s)} ds + c \nh{v_{0, n}} + c\nh{u_{0, n}} \, .
\]
We deal with $T_7$ by direct integration and the use of the H\"older inequality:
\[
T_7(t) = \alpha\scal{\nabla w_n(t)}{\nabla w_{0, n}}_H - \alpha\nh{\nabla w_{0, n}} \leq \frac{\alpha}{8} \nh{\nabla w_n(t)} + \alpha\nh{\nabla w_{0, n}} \, .
\]
Now we pay attention to $T_4$ and integrate by parts in time:{\cred
\[
\begin{split}
T_4(t) = \dual{1*f_n^{(1)}(t)}{w_n(t)} - \int_0^t \dual{f_n^{(1)}(s)}{w_n(s)}ds 
\leq \frac{1}{2\sigma} \norm{1*f_n^{(1)}(t)}^2_{V'}  \\ + \frac{\sigma}{2} \norm{w_n(t)}^2_V
+ \frac{1}{2}\int_0^t\norm{f_n^{(1)}(s)}^2_{V'}ds + \frac{1}{2}\int_0^t \norm{w_n(s)}^2_V ds \, ,
\end{split}
\]
}%
where $\sigma> 0$ is arbitrary, to be set later. According to the definition of the norm in $V$ and the inequality \eqref{fond. t. calculus}, we have
\[
\begin{split}
T_4(t)\leq \frac{1}{2\sigma} \norm{1*f_n^{(1)}(t)}^2_{V'} + \sigma T\int_0^t\nh{\partial_t w_n(s)}ds + \frac{\sigma}{2}\nh{\nabla w_n(t)} \\
+ \frac{1}{2}\int_0^t\norm{f_n^{(1)}(s)}^2_{V'}ds + T\int_0^t\left(\int_0^s \nh{\partial_t w_n(\tau)}d\tau\right)ds \\
+ \frac12 \int_0^t\nh{\nabla w_n(s)}ds + T\left(\sigma + 1\right)\nh{w_{0, n}} \, .
\end{split}
\]
We collect all the terms containing $\norm{\partial_t w_n}_{L^2(0, t; \, H)}$ and $\norm{\nabla w_n(t)}_H$ in the left-hand side of \eqref{test 1, 1}; their coefficients turn out to be, respectively,
\[
k_1 = \frac{1}{2} - T\sigma \, , \quad k_2 = \frac{1}{2}\left(\frac{\alpha}{2} - \sigma \right) \, .
\]
We choose $\sigma \ {\cred \leq} \   \min\left\{\alpha/4, \, 1/4T \right\}$, so that $k_1 \geq 1/4$, 
$k_2 \geq \alpha/8$. We also remark that the assumptions \eqref{f convergence} and \eqref{data 
convergence} enable us to get a bound for terms involving ${\cred f_n^{(1)}, \,f_n^{(2)} } $ 
and the initial data. Finally, adding 
\eqref{test 1, 2} and \eqref{test 1, 1} and taking into account all the previous inequalities, we obtain
\[
\begin{split}
\frac{1}{2}\nh{u_n(t)} + \int_0^t \nh{\nabla u_n(s)} ds + \frac{1}{4} \int_0^t \nh{\partial_tw_n(s)} ds + \frac{\alpha}{8}\nh{\nabla w_n(t)} \\
\leq c\int_0^t \nh{u_n(s)} ds + T\int_0^t\left(\int_0^s \nh{\partial_t w_n(\tau)}d\tau\right)ds + c_\alpha \int_0^t\nh{\nabla w_n(s)}ds + c_\alpha \, .
\end{split}
\]
The Gronwall lemma entails{\cred
\begin{equation}
\norm{u_n}_{\vett{L^\infty}{H} \cap \vett{L^2}{V}} +
\norm{w_n}_{\vett{H^1}{H}} + \sqrt{\alpha} \norm{w_n}_{\vett{L^\infty}{V}} \leq c_{\alpha} \, .
\label{estimate-1}
\end{equation}
}%

\textbf{Second a priori estimate.} Since $\phi_\varepsilon$ is at most of quadratic growth, by definition, and $\gamma_\varepsilon$ is Lipschitz-continuous, from the estimate \eqref{estimate-1} we directly derive
\begin{equation}
\norm{\phi_\varepsilon(u_n)}_{\vett{L^\infty}{L^1(\Omega)}} \leq c'_{\alpha, 1}
\label{estimate 2, 1}
\end{equation}
\begin{equation}
\norm{\gamma_\varepsilon(u_n)}_{L^2(Q)} \leq c'_{\alpha, 2}\, ;
\label{estimate 2, 2}
\end{equation}
where the symbols $c'_{\alpha, i}$ denote positive constants, possibly depending on $\varepsilon$ and 
$\alpha$, but not on  $n$ and $\beta$.

{\cblu By \eqref{numeroform}, we can easily check that
\[
\scal{y}{z}_H = \scal{P_ny}{z}_H \qquad \textrm{for all } y\in V\, , \quad z\in V_n 
\]
where $P_ny$ is the projection of $y$ in $V_n$, with respect to $V$. Then,}
as we have a uniform estimate for $u_n$ in $\vett{L^2}{V}$, it is not difficult to extract from \eqref{eq. B n} the property
\begin{equation}
\norm{\partial_t u_n}_{\vett{L^2}{V'}} \leq c'_{\alpha, 3} \, .
\label{estimate 2, 3}
\end{equation}

\textbf{Third a priori estimate.} We take $v = \partial_t w_n(t)$ as a test function in equation \eqref{eq. A n} and integrate over $(0, t)$; {\cred thanks to} the H\"older inequality, we get
\begin{equation}
\begin{split}
\frac{1}{2}\nh{\partial_t w_n(t)} + \alpha\int_0^t\nh{\nabla\partial_t w_n(s)}ds + \frac{\beta}{2} \nh{\nabla w_n(t)} \hskip2.5cm \\
\leq \int_0^t \dual{f_n^{(1)} - \partial_t u_n(s)}{\partial_t w_n(s)} ds  
+ \int_0^t\norm{f_n^{(2)}(s)}_H \norm{\partial_t w_n(s)}_H ds \\
+ \frac{1}{2}\nh{v_{0,\, n}} + \frac{\beta}{2}\nh{\nabla w_{0, n}} \, .
\end{split}
\label{test 3, 1}
\end{equation}
We consider the term involving $f_n^{(1)} - \partial_t u_n$:
\[
\begin{split}
\int_0^t \dual{f_n^{(1)} - \partial_t u_n(s)}{\partial_t w_n(s)} ds \leq \frac{c}{\alpha} \int_0^t\norm{f_n^{(1)}(s)}^2_{V'}ds + \frac{c}{\alpha} \int_0^t\norm{\partial_t u_n(s)}^2_{V'}ds \\
+ \frac{\alpha}{2} \int_0^t\nh{\partial_t w_n(s)} ds + \frac{\alpha}{2} \int_0^t \nh{\nabla\partial_t w_n(s)} ds \, .
\end{split}
\]
Because of the estimate \eqref{estimate 2, 3} {\cred and the properties} \eqref{f convergence} and \eqref{data convergence}, from \eqref{test 3, 1} we deduce 
\[
\begin{split}
\frac{1}{2}\nh{\partial_t w_n(t)} + \frac{\alpha}{2} \int_0^t\nh{\nabla\partial_t w_n(s)}ds + \frac{\beta}{2} \nh{\nabla w_n(t)} \\
\leq c' + \frac{\alpha}{2} \int_0^t\nh{\partial_t w_n(s)}ds + \int_0^t\norm{f_n^{(2)}(s)}_H \norm{\partial_t w_n(s)}_H ds \, ,
\end{split}
\]
where $c'$ depends on $\varepsilon, \, \alpha$. Hence, by a generalized version of the Gronwall
lemma {\cred (see, e.g., \cite[pp. 156--157]{Brezis})}, we {\cred infer that}
\begin{equation}
\norm{w_n}_{\vett{W^{1, \, \infty}}{H}} + {\cred \sqrt{\alpha}} \norm{w_n}_{\vett{H^1}{V}} \leq c'_{\alpha, 4} \, .
\label{estimate 3, 1}
\end{equation}

%
%

\textbf{Passage to the limit as $n \longrightarrow + \infty$.} From the estimates \eqref{estimate-1}, \eqref{estimate 2, 1}--\eqref{estimate 2, 3}, \eqref{estimate 3, 1}, with standard arguments of weak or weak* compactness we can find functions $(w_\varepsilon, \, u_\varepsilon)$ such that, possibly taking a subsequence as $n\longrightarrow + \infty$, {\cred
\begin{eqnarray}
 w_n \rightharpoonup^* w_\varepsilon && \textrm{in } \ \vett{W^{1, \, \infty}}{H} \cap \vett{L^\infty}{V} \label{1conv}\\
w_n \rightharpoonup   w_\varepsilon && \textrm{in }\  \vett{H^1}{V} \label{2conv}\\
u_n \rightharpoonup   u_\varepsilon && \textrm{in } \ \vett{H^1}{V'} \cap \vett{L^2}{V} \label{3conv}\\
u_n \rightharpoonup^*   u_\varepsilon && \textrm{in } \ \vett{L^\infty}{H}  \, .\label{4conv}
\end{eqnarray}
Note that \eqref{2conv} implies 
}%
the strong convergence 
\begin{equation}
w_n \longrightarrow w_\varepsilon \qquad \textrm{in } C^0\left([0, T]; \, H\right) ;
\label{w C0}
\end{equation}
on the other hand, {\cred the generalised Ascoli theorem and the Aubin-Lions lemma 
(see, e.g., \cite[pp.~57--58]{Lions} and  \cite[Sect.~8, Cor.~4]{Simon})
entail
\begin{equation}
u_n \longrightarrow u_\varepsilon \quad \textrm{ strongly in } C^0\left([0, T]; \, V' \right)  
\textrm{ and in } L^2(Q) ;
\label{u C0}
\end{equation}
thus}, since $g$ and $\gamma_\varepsilon$ are Lipschitz-continuous, 
we easily check that
\[
g(u_n) \longrightarrow g(u_\varepsilon) \quad \textrm{ {\cred and} } \quad  \gamma_\varepsilon(u_n) \longrightarrow \xi_\varepsilon \quad {\cred \textrm{ strongly in } L^2(Q),} 
\]
{\cred where $\xi_\varepsilon = \gamma_\varepsilon(u_\varepsilon )$. We then take the limit as 
$n\longrightarrow + \infty$ in \eqref{eq. A n}--\eqref{in.cond}}
and see that $(w_\varepsilon, \, u_\varepsilon, \, \xi_\varepsilon)$ fulfill{\cred s} equations 
{\cred \eqref{xi}--\eqref{initial condition}, where $\gamma $ is replaced by $\gamma_\varepsilon$.
Indeed, 
by \eqref{w C0}--\eqref{u C0} and \eqref{data convergence}}, it is obvious that $w_\varepsilon(0) = w_0$, $u_\varepsilon(0) = 
u_0$.
To deal with the last initial condition properly, we fix a test function $v\in V_m$, where $m\geq 1$ is arbitrary, and we integrate in time equation \eqref{eq. A n}; we get equation \eqref{new A n}, {\cred for $0 \leq t \leq T$ and} $n\geq m$. Arguing as in {\cred \cite[pp.~12--13]{Lions}}, we can take the limit in \eqref{new A n}, \eqref{eq. B n} and check that $(w_\varepsilon, \, u_\varepsilon, \, \xi_\varepsilon)$ fulfills
\begin{equation}
\begin{split}
\dual{\partial_tw_\varepsilon(t)}{v} = - \alpha\scal{\nabla w_\varepsilon(t)}{\nabla v}_H - \beta\scal{1*\nabla w_\varepsilon(t)}{\nabla v}_H \\
- \dual{u_\varepsilon(t)}{v} + \dual{1*f(t)}{v} + \alpha\scal{\nabla w_0}{\nabla v}_H + \scal{v_0 + u_0}{v}_H
\end{split}
\label{eq. A int}
\end{equation}
\begin{equation}
\dual{\partial_tu_\varepsilon(t)}{v} + \scal{\nabla u_\varepsilon(t)}{\nabla v}_H + \scal{\xi_\varepsilon(t)}{v}_H + \scal{g(u_\varepsilon)(t)}{v}_H 
= \scal{\partial_tw_\varepsilon(t)}{v}_H
\label{eq. B int}
\end{equation}
for {\cred a.a. $t\in(0,T)$,} $m\geq 1$ and $v\in V_m$; by a density argument, the same equalities hold when $v\in V$. Since the right-hand side in \eqref{eq. A int} is a continuous function in $[0, T]$, taking $t = 0$ we find that
\[
{\cred \dual{\partial_t w_\varepsilon(0)}{v}} = \scal{v_0}{v}_H \qquad \textrm{for all } v\in V 
\]
{\cred whence the second of \eqref{initial condition} follows.}

\bigskip
\textbf{Fifth a priori estimate.} As a consequence of the weak lower semi-con\-ti\-nui\-ty of the norm in a Banach space, $(w_\varepsilon, \, u_\varepsilon, \, \xi_\varepsilon)$ satisfy the estimate \eqref{estimate-1}; we now need to improve estimates \eqref{estimate 2, 1}--\eqref{estimate 2, 3}, \eqref{estimate 3, 1}.
 
We first notice that, because of the Lipschitz-continuity of $\gamma_\varepsilon$, $\xi_\varepsilon(t)\in V$ for all $t$; thus, we can choose $v = \xi_\varepsilon(t)$ in equation \eqref{eq. B int} and integrate over $(0, t)$, to get
\begin{equation}
\begin{split}
\int_{Q_t} \partial_t u_\varepsilon \,\xi_\varepsilon + \int_{Q_t} \gamma_\varepsilon'(u_\varepsilon)\abs{\nabla u_\varepsilon}^2 + \int_0^t\nh{\xi_\varepsilon(s)}ds = \int_{Q_t} g(u_\varepsilon)\,\xi_\varepsilon + \int_{Q_t} \partial_t w_\varepsilon \,\xi_\varepsilon \, .
\end{split}
\label{test ?}
\end{equation}
{\cred In view of \eqref{prope},} we have
\[
\int_{Q_t} \partial_t u_\varepsilon \, \xi_\varepsilon = \int_{Q_t} \frac{\partial}{\partial t}\left(\phi_\varepsilon(u_\varepsilon)\right) = \norm{\phi_\varepsilon(u_\varepsilon(t))}_{L^1(\Omega)} - \norm{\phi_\varepsilon(u_0)}_{L^1(\Omega)} \, ;
\]
on the other hand, because of the Lipschitz continuity of $g$,
\[
\int_{Q_t} g(u_\varepsilon)\xi_\varepsilon \leq c \int_{Q_t} \left(\abs{u_\varepsilon} + 1\right)\xi_\varepsilon \leq c\int_0^t \left(\nh{u_\varepsilon(s)} + 1\right)ds + \frac{1}{2}\int_0^t \nh{\xi_\varepsilon(s)} ds \, .
\]
From these estimates and \eqref{test ?}, we derive{\cred
\[
\begin{split}
\int_\Omega \phi_\varepsilon(u_\varepsilon)(t) + \int_{Q_t}\gamma_\varepsilon'(u_\varepsilon)\abs{\nabla u_\varepsilon}^2 + 
\frac{1}{2}\int_0^t\nh{\xi_\varepsilon(s)}ds \hskip1.5cm\\
{}\leq c\int_0^t\nh{u_\varepsilon(s)}ds + c\int_0^t\nh{\partial_t w_\varepsilon(s)}ds +  \int_\Omega\phi_\varepsilon(u_0) + c \, .
\end{split}
\]
}%
We notice that the second term in the lef-hand side is nonnegative, because of the monotonicity of $\gamma_\varepsilon$. Secondly, accounting for \eqref{estimate-1}, 
{\cred \eqref{prope} and} \eqref{initial data}, {\cred we infer that 
\begin{equation}
\norm{\phi_\varepsilon(u_\varepsilon)}_{\vett{L^\infty}{L^1(\Omega){1}}} +
\norm{\gamma_\varepsilon(u_\varepsilon)}_{L^2(Q)} \leq c_\alpha \, .
\label{estimate-4}
\end{equation}
}%
Now, by comparison in the equation \eqref{eq. B int}, we have
\begin{equation}
\norm{\partial_t u_\varepsilon}_{\vett{L^2}{V'}} \leq {\cred c_\alpha} \, ;
\label{estimate 4, 3}
\end{equation}
{\cred and consequently
we can also establish the estimate \eqref{estimate 3, 1}, now} for a constant which is independent of $\varepsilon$.

\bigskip
\textbf{Passage to the limit as $\varepsilon \searrow 0$.} We are able to repeat the compactness argument {\cred as above} and find $(w, \, u, \, \xi)$, {\cred a candidate for the solution to} Problem~$\Pab$, as a limit of a subsequence of $(w_\varepsilon, \, u_\varepsilon, \, \xi_\varepsilon)$. The proof will be easily completed by {\cred the} passage to the limit as $\varepsilon \searrow 0$, provided that we deduce \eqref{xi}.

By construction, we can assume that
\[
\xi_\varepsilon \rightharpoonup \xi \ \textrm{ in } L^2(Q) \, , \qquad u_\varepsilon \longrightarrow u \ \textrm{ in } L^2(Q) \, , 
\]
from which the equality
\[
\lim_{\varepsilon\searrow 0}\int_Q \xi_\varepsilon u_\varepsilon = \int_Q \xi u 
\]
follows; {\cred at this point, we apply \cite[Prop.~1.1, p.~42]{Barbu} and deduce \eqref{xi}. Thus,} the proof of the existence of a solution {\cred to} Problem~$\Pab$ is complete.

\section{Regularity and strong solutions}
\label{reg1}
This section is devoted to the derivation of further a priori estimates on the {\cred approximating} solutions $(w_n, \, u_n, \, \xi_n)$, which are independent of $n$ and $\varepsilon$, under stronger assumptions. The same compactness -- passage to the limit arguments then apply, and this will prove Theorem \ref{th: strong solution}. We first notice that the hypothesis {\cred \eqref{initial data strong} and  $V_n \subseteq W$} make 
it possible to assume
\begin{equation}
{\cred w_{0, n} \longrightarrow w_0 \  \textrm{ in } W\, , \qquad v_{0, n} 
\longrightarrow v_0  \ \hbox{ and } \  u_{0, n} \longrightarrow u_0 \ \textrm{ in } V \, ;}
\label{piuconvstr}
\end{equation}
on the other hand, owing to \eqref{f strong}, we can require $f_n^{(1)}\in L^2(Q)$, $f_n^{(2)}\in \vett{L^1}{V}$ for all $n\in\N$ and
\begin{equation}
f_n^{(1)} \longrightarrow f^{(1)} \  \textrm{ in } L^2(Q)\, , \qquad
f_n^{(2)} \longrightarrow f^{(2)} \  \textrm{ in } \vett{L^1}{V} .
\label{fpiuconvstr}
\end{equation}

\bigskip
\textbf{Sixth a priori estimate.} We choose $v = \partial_t w_n(t)$ in the equation \eqref{eq. A n} and integrate over $(0, t)$; an application of the H\"older inequality yields
\begin{equation}
\begin{split}
\frac{1}{2}\nh{\partial_t w_n(t)} + \alpha\int_0^t\nh{\nabla\partial_t w_n(s)}ds + \frac{\beta}{2}\nh{\nabla w_n(t)} \leq - \int_{Q_t}\partial_t u_n\partial_t w_n \\
+ \int_0^t\norm{f_n(s)}_H \norm{\partial_t w_n(s)}_H ds + \frac{1}{2}\nh{v_{0, n}} + \frac{\beta}{2}\nh{\nabla w_{0, n}} \, .
\end{split}
\label{test 5, 1}
\end{equation}
Now{\cred , we take} $v = \partial_t u_n(t)$ in \eqref{eq. B n} and integrate over $(0, t)$; recalling that $\gamma_\varepsilon = \phi_\varepsilon'$, using the H\"older inequality and the Lipschitz-continuity of $g$, we get
\begin{equation}
\begin{split}
\frac{1}{2} \int_0^t\nh{\partial_t u_n(s)}ds + \frac{1}{2}\nh{\nabla u_n(t)} + \norm{\phi_\varepsilon(u_n(t))}_{L^1(\Omega)} \\
\leq \int_{Q_t} \partial_t u_n \, \partial_t w_n + c\int_0^t \left(\nh{u_n(s)} + 1\right)ds + \norm{\phi_\varepsilon(u_{0, n})}_{L^1(\Omega)} \, .
\end{split}
\label{test 5, 3}
\end{equation}
Adding \eqref{test 5, 1} and \eqref{test 5, 3}, thanks to the assumptions \eqref{initial data}, \eqref{data convergence}, the inequality \eqref{fond. t. calculus} and $\phi_\varepsilon \leq \phi$, we finally have
\[
\begin{split}
\frac{1}{2}\nh{\partial_t w_n(t)} + \alpha\int_0^t\nh{\nabla\partial_t w_n(s)}ds + \frac{\beta}{2}\nh{\nabla w_n(t)} \\
+ \frac{1}{2} \int_0^t\nh{\partial_t u_n(s)}ds + \frac{1}{2}\nh{\nabla u_n(t)} + \norm{\phi_\varepsilon(u_n(t))}_{L^1(\Omega)} \\
\leq c\int_0^t \left(\int_0^s\nh{\partial_t u_n(\tau)}d\tau\right)ds + \int_0^t \norm{f_n(s)}_H\norm{\partial_t w_n (s)}_H ds + c \, .
\end{split}
\]
The generalised Gronwall lemma {\cred (see, e.g., \cite[pp. 156--157]{Brezis})} enable{\cred s} us to achieve
\begin{equation}
\norm{w_n}_{\vett{W^{1, \infty}}{H}} + {\cred \sqrt{\alpha}} \norm{w_n}_{\vett{H^1}{V}} + {\cred \sqrt{\beta}} \norm{w_n}_{\vett{L^\infty}{V}} \leq c_1
\label{estimate 5, 1}
\end{equation}
\begin{equation}
\norm{u_n}_{\vett{H^1}{H} \cap \vett{L^\infty}{V}} \leq c_2.
\label{estimate 5, 2}
\end{equation}

\begin{remark} Only the hypotheses \eqref{alpha, beta}--\eqref{initial data} and $f\in\vett{L^1}{H}$ have been effectively exploited in the proof of this estimate.
\end{remark}
\begin{remark} By means of \eqref{estimate 5, 1}--\eqref{estimate 5, 2}, the estimates \eqref{estimate-4}--\eqref{estimate 4, 3} can be {\cred rewritten} in terms of some constant which is independent of $\alpha$.
\end{remark}

\bigskip
\textbf{Seventh a priori estimate.} We take $v = -\Delta u_n(t)$ in equation \eqref{eq. B n}; this is possible, because of the special choice of the approximating space $V_n$. We integrate over $(0, t)$ and use the H\"older inequality and the Lipschitz continuity of $g$:
\[
\begin{split}
\frac{1}{2}\nh{\nabla u_n(t)} + \int_0^t \nh{\Delta u_n(s)} ds + \int_{Q_t}\gamma'_\varepsilon(u_n)\abs{\nabla u_n}^2 \\
= - \int_{Q_t} g'(u_n)\abs{\nabla u_n}^2 - \int_{Q_t} \partial_t w_n \Delta u_n + \frac{1}{2}\nh{\nabla u_{0, n}} \\
\leq c\norm{\nabla u_n}^2_{\vett{L^2}{H}} + \frac{1}{2}\norm{\partial_t w_n}^2_{\vett{L^2}{H}} + \frac{1}{2}\int_0^t\nh{\Delta u_n(s)} ds + \frac{1}{2}\nh{\nabla u_{0, n}} \, .
\end{split}
\]
The monotonicity of $\gamma_\varepsilon$ yields that the last term in the lef-hand side is non negative. Owing to {\cred conditions \eqref{piuconvstr}} on the data and estimates \eqref{estimate 5, 1}, \eqref{estimate 5, 2}, we have
\[
\frac{1}{2}\int_0^t\nh{\Delta u_n(s)} ds \leq c \qquad \textrm{for all } 0\leq t \leq T \, ;
\]
hence, on account of this inequality, the estimate \eqref{estimate 5, 2} and the boundary conditions for $u_n$, {\cred known} regularity results for elliptic problems entail
\begin{equation}
\norm{u_n}_{\vett{L^2}{W}} \leq c_3 \, ,
\label{estimate 6, 1}
\end{equation}
where ${\cred c_3}$ does not depend on $\alpha$, $\beta$.

\bigskip
\textbf{Eigth a priori estimate.} Since $w_n \in C^2([0, T]; \, V_n)$, the special choice of $V_n$ enable{\cred s} us to take $v = - \Delta \partial_t w_n(t)$ as a test function in the equation \eqref{eq. A n}. We integrate over $(0, t)$ and use the H\"older inequality:
\begin{equation}
\begin{split}
\frac{1}{2}\nh{\nabla\partial_t w_n(t)} + \alpha\int_0^t\nh{\Delta\partial_t w_n(s)}ds + \frac{\beta}{2}\nh{\Delta w_n(t)} \\
{}\leq {\cred  \frac{\alpha}{2}\int_0^t\nh{\Delta\partial_t w_n(s)} ds + \frac{1}{\alpha} \int_0^t \nh{\partial_t u_n(s)} ds + \frac{1}{\alpha}\int_0^t\nh{f^{(1)}_n(s)}ds }\\
{}- \int_{Q_t} f_n^{(2)}\Delta\partial_t w_n + \frac{1}{2}\nh{\nabla v_{0, n}} + \frac{\beta}{2}\nh{\Delta w_{0, n}} \, .
\end{split}
\label{test 7, 1}
\end{equation}
For the term involving $f_n^{(2)}$, we integrate by parts in space, {\cred recalling} that $\partial_n v = 0$ for all $v\in V_n$:
\begin{equation}
\abs{\int_{Q_t} f_n^{(2)}\Delta\partial_t w_n} = \abs{\int_{Q_t}\nabla f_n^{(2)}\cdot\nabla 
\partial_t w_n} \leq \int_0^t \norm{\nabla f_n^{(2)}{\cred (s)}}_H \norm{\nabla \partial_t 
w_n(s)}_H ds \, .
\label{test 7,1bis}
\end{equation}
{\cred Then, in view of \eqref{piuconvstr}, \eqref{fpiuconvstr}, \eqref{estimate 5, 2} and 
owing to the generalized Gronwall lemma (see \cite[pp. 156--157]{Brezis}), from 
\eqref{test 7, 1}--\eqref{test 7,1bis} we obtain}
\begin{equation}
\norm{w_n}_{\vett{W^{1, \infty}}{V}} + {\cred \sqrt{\alpha}} \norm{w_n}_{\vett{H^1}{W}} + {\cred \sqrt{\beta}} \norm{w_n}_{\vett{L^\infty}{W}} \leq c_{\alpha, 4} \, .
\label{estimate 7, 1}
\end{equation}
Finally, if we choose $v = \partial_t^2 w_n(t)$ in the equation \eqref{eq. A n}, we get
\[
\cred{ \nh{\partial_t^2 w_n(t)} \leq \left\{ \alpha\norm{\partial_t w_n(t)}_W + \beta\norm{w_n(t)}_W + 
\norm{\partial_t u_n(t)}_H + \norm{f_n(t)}_H \right\} \norm{\partial_t^2 w_n(t)}_H \, ;}
\]
thanks to the estimates above, it is easy to derive
\begin{equation}
\norm{\partial_t^2 w_n}_{\vett{L^1}{H}} \leq {\cred c_{\alpha, 5} } \, .
\label{estimate 7, 3}
\end{equation}
Having established all the a priori estimates corresponding to \eqref{w strong}--\eqref{u strong} on the solutions of the approximating problem, we have completed the proof of Theorem \ref{th: strong solution}.

\section{Further regularity}
\label{reg2}
Throughout this section we assume \eqref{u0 strong} {\cred and} \eqref{u0 gamma} in addition to all the hypotheses we had in {\cred Section~\ref{reg1}}. As we are interested in {\cred proving} Theorem~\ref{th: regularity}, we should get further estimates on the solution of the approximated problem. By the stronger assumptions on the initial data, we can require
\begin{equation}
u_{0, n} \longrightarrow u_0 \quad \textrm{ in } W \, .
\label{u0 convergence strong}
\end{equation}

Consider the equation \eqref{eq. B n} and derive it, with respect to time, obtaining
\[
\begin{split}
\scal{\partial_t^2 u_n(t)}{v}_H + \scal{\nabla\partial_t u_n(t)}{\nabla v}_H + \scal{\gamma'_\varepsilon(u_n(t))\partial_t u_n(t)}{v}_H \\
+ \scal{g'(u_n(t))\partial_t u_n(t)}{v}_H = \scal{\partial_t^2 w_n(t)}{v}_H
\end{split}
\]
for all $v\in V_n$ and a.a. $\cred t\in (0,T)$. We choose $v = \partial_t u_n(t)$ as an admissible test function, integrate over $(0, t)$ and use the Lipschitz continuity of $g$ to get
\begin{equation}
\begin{split}
\frac{1}{2}\nh{\partial_t u_n(t)} + \int_0^t \nh{\nabla \partial_t u_n(s)} ds + \int_{Q_t} \gamma'_\varepsilon(u_n)\abs{\partial_t u_n}^2 {\cred {}\leq c\int_0^t \nh{\partial_t u_n(s)} ds}   \\
{}+ \int_0^t \norm{\partial_t^2 w_n(s)}_H\norm{\partial_t u_n(s)}_H ds + \frac{1}{2} \nh{\partial_t u_n(0)} \, .
\end{split}
\label{test 8, 1}
\end{equation}
Since the last term in the left-hand side is non negative because of {\cred the 
monotonicity of $\gamma_\varepsilon$}, 
if we had a bound for the last term in the right-hand side, we could use the generalized 
Gronwall lemma to conclude. In order to provide such an estimate, we set $t = 0$, $v = 
\partial_t u_n(0)$ in the equation \eqref{eq. B n}; we obtain
\[
\nh{\partial_t u_n(0)} \leq \left\{ \norm{\Delta u_{0, n}}_H + 
\norm{\gamma_\varepsilon(u_{0, n})}_H + \norm{g(u_{0, n})}_H + \norm{v_{0, n}}_H 
\right\} \norm{\partial_t u_n(0)}_H
\]
and thus, taking into account the Lipschitz {\cred continuity of $g$, we infer}
\[
\begin{split}
\norm{\partial_t u_n(0)} \leq \norm{\Delta u_{0, n}}_H + 
\norm{\gamma'_\varepsilon}_{L^\infty(\R)}\norm{u_{0, n} - u_0}_H + 
\norm{\gamma_\varepsilon(u_0)}_H\\
 + c \left( \norm{u_{0, n}}_H + 1\right) + \norm{v_{0, n}}_H .
\end{split}
\]
Now, assumptions \eqref{u0 convergence strong} and \eqref{data convergence}, as well as \eqref{u0 gamma} and $\abs{\gamma_\varepsilon} \leq \abs{\gamma^0}$, enable us to achieve
\begin{equation}
\norm{\partial_t u_n(0)} \leq c
\label{test 8, 2}
\end{equation}
for all $\varepsilon > 0$ and $n$ large enough, depending on $\varepsilon$; these requests on parameters are not restrictive, as we first take the limit for $n \longrightarrow + \infty$, then for $\varepsilon \searrow 0$. 
From \eqref{test 8, 1} and \eqref{test 8, 2} we {\cred deduce that}
\begin{equation}
\norm{u_n}_{\vett{W^{1, \infty}}{H} \cap \vett{H^1}{V}} \leq {\cred c_{\alpha, 6}} 
\, .
\label{estimate 8, 1}
\end{equation} 

Finally, we consider equation \eqref{eq. B n} and we rewrite it in the form
\[
\scal{\nabla u_n(t)}{\nabla v}_H + \scal{\gamma_\varepsilon(u_n(t))}{v}_H = \scal{F_n(t)}{v} \, , 
\] 
for all $v\in V_n$ and a.a. {\cred $t\in (0,T)$, where $F_n = \partial_t w_n - \partial_t u_n - g(u_n) $}. Testing with $v = -\Delta u_n(t)$ the previous equation and integrating by parts in space, we obtain
\[
\nh{\Delta u_n(t)} + \int_\Omega \gamma'_\varepsilon(u_n(t))\abs{\nabla u_n(t)}^2 \leq \norm{F_n(t)}_H \norm{\Delta u_n(t)}_H \qquad \textrm{for all } 0\leq t \leq T \, .
\]
Since the estimates \eqref{estimate 5, 1} and \eqref{estimate 8, 1} entail
\[
\norm{F_n}_{\vett{L^\infty}{H}} \leq c
\]
and we can apply the regularity results for elliptic problems, we deduce
\begin{equation}
\norm{u_n}_{\vett{L^\infty}{W}} \leq {\cred c_{\alpha, 7}} \, ,
\label{estimate 8, 2}
\end{equation}
thus concluding the proof of Theorem \ref{th: regularity}.

\section{$L^\infty$ estimates}
\label{reg3}
The aim of this section is to obtain $L^\infty$ estimates on $w_t$ and on $\xi$, under the hypotheses \eqref{f stronger} and \eqref{gamma strong}.

We first deal with $w_t$. Setting $\varphi = \alpha w_t + \beta w$, Theorem \ref{th: strong solution} entail that the equalities
\[
\frac{1}{\alpha}\varphi_t - {\cred \Delta\varphi } = \frac{\beta}{\alpha} w_t - u_t + f \quad \textrm{ in } Q \, , \qquad 
\partial_n \varphi = 0 \quad \textrm{ on }\Gamma \times (0, T)
\]
hold almost everywhere. Furthermore, the assumption \eqref{f stronger}, the estimates \eqref{estimate 5, 1} and {\cred \eqref{estimate 8, 1}} and the continuous embedding $V \hookrightarrow L^6(\Omega)$  {\cred (valid if  $\Omega\subseteq\R^3$ is a bounded Lipschitz domain)}, yield
\[
\frac{\beta}{\alpha} w_t - u_t + f \in \vett{L^\infty}{H} + \vett{L^r}{L^6(\Omega)} \, , \qquad \textrm{with } r > 4/3 \, .
\]
In these conditions, {\cred Theorem~7.1 in \cite[p.~181]{LSU}} applies and {\cred ensures that} $\varphi\in L^\infty(Q)$. Since we already know {\cred that
$w\in L^\infty(Q)$ (as it is implied, for example, by \eqref{estimate 7, 1})}, we have $w_t\in L^\infty(Q)$ and
\begin{equation}
\norm{w_t}_{L^\infty(Q)} \leq \frac{1}{\alpha}\norm{\varphi}_{L^\infty(Q)} + \frac{c\beta}{\alpha} \norm{w}_{\vett{{\cred L^\infty}}{W}} \leq c_{\alpha, 8} \, .
\label{estimate 9, 1}
\end{equation}
We notice that, being $\alpha$ fixed and letting $\beta$ vary in a bounded set, we can find an upper bound for the constant $c_{\alpha, 8}$.

\bigskip
In order to prove a $L^\infty$ estimate for $\xi$, we consider the solution $(w_\varepsilon, \, u_\varepsilon)$ to the approximating problem, in which the Yosida regularization appears; we then fix $p\in(1, \, +\infty)$ and get a bound for $\norm{\gamma_\varepsilon(u_\varepsilon)}_{L^p(Q)}$, which is independent of $p$, $\varepsilon$. From this, we will obtain a uniform bound for 
\[
\norm{\xi_\varepsilon}_{L^\infty(Q)} = \lim_{p\rightarrow+\infty} \norm{\gamma_\varepsilon(u_\varepsilon)}_{L^p(Q)} \, ,
\]
and, via a weak* compactness argument, $\xi\in L^\infty(Q)$. For the sake of simplicity, we do not plug in the subscript $\varepsilon$ in the solution any more.

We know that the equalities
\begin{equation}
u_t - \Delta u + \gamma_\varepsilon(u) + g(u) = w_t \quad \textrm{in }Q \, , 
\label{u_vareps}
\end{equation}
\[
\partial_n u = 0 \quad \textrm{on }\Gamma\times (0, T) \, , \qquad u(0) = u_0 \quad \textrm{in } \Omega
\]
hold a.e.; we choose $\abs{\gamma_\varepsilon(u)}^{p - 1}\gamma_\varepsilon(u)$ as a test function, by which we multiply both sides of the equation \eqref{u_vareps} -- this is admissible since $u\in L^\infty(Q)$. Integrating over $Q$, we get
\begin{equation}
\begin{split}
\int_Q \frac{\partial}{\partial t}\phi_{\varepsilon, \, p}(u) + \int_Q \nabla u \cdot \nabla\left(\abs{\gamma_\varepsilon(u)}^{p - 1}\gamma_\varepsilon(u)\right) + \int_Q \abs{\gamma_\varepsilon(u)}^{p + 1} \\
= \int_Q \left(w_t - g(u)\right)\abs{\gamma_\varepsilon(u)}^{p - 1}\gamma_\varepsilon(u) \, ,
\end{split}
\label{test 9, 1}
\end{equation}
where we have set
\[
\phi_{\varepsilon, \, p}(t) = \int_0^t \abs{\gamma_\varepsilon(s)}^{p - 1}\gamma_\varepsilon(s)\, ds \qquad \textrm{for all } t\in\R \, ;
\]
$\gamma_\varepsilon$ is increasing and $\gamma_\varepsilon(0) = 0$, so we have $\phi_{\varepsilon, \, p}\geq 0$ for all $\varepsilon$, $p$. Since $w_t, \, u\in L^\infty(Q)$ and $g$ is continuous, for the right-hand side we have
\[
\abs{\int_Q \left(w_t - g(u)\right)\abs{\gamma_\varepsilon(u)}^{p - 1}\gamma_\varepsilon(u)} \leq c_\alpha \norm{\gamma_\varepsilon(u)}^p_{L^p(\Omega)} \, ;
\]
on the other hand, a direct calculation and the monotonicity of $\gamma_\varepsilon$ show that
\[
\nabla u \cdot \nabla\left(\abs{\gamma_\varepsilon(u)}^{p - 1}\gamma_\varepsilon(u)\right) = p\gamma'_\varepsilon(u)\abs{\gamma_\varepsilon(u)}^{p - 1}\abs{\nabla u}^2 \geq 0  \qquad \textrm{a.e. in } Q \, . 
\]
Collecting all the information we have {\cred obtained} so far, from \eqref{test 9, 1} we derive
\begin{equation}
\int_\Omega \phi_{\varepsilon, \, p}(u(T)) + \norm{\gamma_\varepsilon(u)}^{p + 1}_{L^{p+1}(Q)} \leq c_\alpha \norm{\gamma_\varepsilon(u)}^p_{L^p(Q)} + \int_\Omega \phi_{\varepsilon, \, p}(u_0)
\label{test 9, 2}
\end{equation}
and, since the first term can be ignored, we need only to find an estimate for the last term. We recall that, for the Yosida approximation of a maximal monotone graph, the inequality
\[
\abs{\gamma_\varepsilon(s)} \leq \abs{\gamma^0(s)} \qquad \textrm{for all } s\in D(\gamma)\, , \quad \varepsilon > 0
\]
holds {\cred (see, e.g., \cite[Prop.~2.6, p.~28]{Brezis})}; according to that, we have 
\[
\begin{split}
\int_\Omega \phi_{\varepsilon, \, p}(u_0) \leq \int_\Omega \abs{\gamma^0(u_0)}^p \abs{u_0}
\leq \frac{p}{p + 1}\int_\Omega \abs{\gamma^0(u_0)}^{p + 1} + \frac{1}{p + 1}\int_\Omega\abs{u_0}^{p + 1} \\
\leq \frac{p}{p + 1}\int_\Omega \abs{\gamma^0(u_0)}^{p + 1} + \frac{1}{p + 1}\norm{u_0}_{L^{p + 1}(\Omega)}^{p + 1}\, , 
\end{split}
\]
where the H\"older and Young inequalities {\cred have been used. We} recall that $u_0\in L^\infty(\Omega)$ by the assumption \eqref{u stronger} {\cred and also}
notice that the same inequalities imply
\[
c_\alpha\norm{\gamma_\varepsilon(u)}^p_{L^p(Q)} \leq \frac{p}{p + 1}\norm{\gamma_\varepsilon(u)}^{p + 1}_{L^{p+1}(Q)} + \frac{c_\alpha}{p + 1} \, .
\]
Now, we come back to the equation \eqref{test 9, 2}; according to the previous estimates, we {\cred infer that}
\[
\frac{1}{p + 1}\norm{\gamma_\varepsilon(u)}^{p + 1}_{L^{p+1}(Q)} \leq \frac{p}{p + 1}\norm{\gamma^0(u_0)}^{p + 1}_{L^{p+1}(\Omega)} + \frac{1}{p + 1}\norm{u_0}_{L^{p + 1}(\Omega)}^{p + 1} + \frac{c_\alpha}{p + 1}
\]
and, hence,
\[
\begin{split}
\norm{\gamma_\varepsilon(u)}_{L^{p+1}(Q)} \leq \left\{p \norm{\gamma^0(u_0)}^{p + 1}_{L^{p+1}(\Omega)} + \norm{u_0}_{L^{p + 1}(\Omega)}^{p + 1} + c_\alpha \right\}^{1/(p + 1)} \\
\leq c_\alpha\left\{ \norm{\gamma^0(u_0)}_{L^\infty(\Omega)} + \norm{u_0}_{L^\infty(\Omega)} + 1\right\} \, ,
\end{split}
\]
which provides the desired estimate and concludes the proof.

\section{Well-posedness of $\Pa$ {\cred and} convergence as $\beta\searrow 0$}
\label{beta=0}
Now we set the notation as in {\cred Section} \ref{sec: Pa}, since we are interested in the proof of Theorems~\ref{th: well-posedness Pa}--\ref{th: second estimate error}. We assume that the hypotheses \eqref{alpha, beta}--\eqref{initial data} are satisfied, and we start by studying the convergence as $\beta\searrow 0$, by a compactness argument.

\bigskip
\textbf{Convergence as $\beta\searrow 0$.} {\cred We recall the a priori estimates \eqref{estimate 3, 1}, \eqref{estimate-1}, \eqref{estimate-4}, \eqref{estimate 4, 3} which are
independent of $\beta$ and thus holding also for  $(w_\beta, \, u_\beta, \, \xi_\beta)$. Moreover, adopting the notation as in \eqref{f splittata}--\eqref{f convergence}, by a comparison in \eqref{eq. A weak} we find out that $\{ \partial_t^2 w_\beta - f^{(2)}_\beta \} $ is uniformly bounded in $\vett{L^2}{V'}$. Therefore, we can find a subsequence $\beta_k\searrow 0$ and functions $w, \, u, \, \xi$ such that
\[
w_{\beta_k} \rightharpoonup^* w \ \hbox{ in } \  \vett{W^{1, \, \infty}}{H} \, , \qquad w_{\beta_k} \rightharpoonup w \ \hbox{ in } \  \vett{H^1}{V}
\]
\[
\partial_t^2 w_\beta - f^{(2)}_\beta \ \rightharpoonup \ w_{tt} - f^{(2)} \ \hbox{ in } \  \vett{L^2}{V'}
\]
\[
u_{\beta_k} \ \hbox{ tends to } \ u \ \hbox{ weakly in } \, \vett{H^1}{V'} \cap \vett{L^2}{V}, \, \hbox{ whence strongly in } \,  L^2(Q), 
\]
\[
\xi_{\beta_k} \rightharpoonup \xi \ \hbox{ in } \  L^2(Q)
\]
as $k\longrightarrow +\infty$, and here part of \eqref{f convergence beta} has been used.
Then, in view of \eqref{g}, \eqref{f convergence beta} and \eqref{data convergence beta}, we can pass to the limit in \eqref{eq. A weak} and \eqref{eq. B weak}, as well as in the initial conditions \eqref{initial condition} which can be recovered weakly in $V'$ at least. On the other hand, $u\in D(\gamma)$ and $\xi\in\gamma(u)$ a.e. in $Q$ follow as a consequence of the above convergences and \cite[Lemma~1.3, p.~42]{Barbu}.}

\bigskip
\textbf{Uniqueness for $\Pa$.} By applying the previous result with $f_\beta = f$, $w_{0, \beta} = w_0$, $v_{0, \beta} = v_0$ and $u_{0, \beta} = u_0$ given, we obtain the existence of a solution to Problem~$\Pa$; we still have to prove the uniqueness. Let $(w_1, \, u_1, \, \xi_1)$ and $(w_2, \, u_2, \, \xi_2)$ be solutions of $\Pa$; we write down the equations for the differences $w = w_1 - w_2$, $u = u_1 - u_2$, $\xi = \xi_1 - \xi_2$ and  integrate with respect to time the first one:
\[
\scal{w_{t}(t)}{v}_H + \alpha\scal{\nabla w(t)}{\nabla v}_H + \scal{u(t)}{v}_H = 0\, , 
\]
\[
\dual{u_t(t)}{v} + \scal{\nabla u(t)}{\nabla v}_H + \scal{\xi(t)}{v}_H + \scal{g(u_1)(t) - g(u_2)(t)}{v}_H = \scal{w_t(t)}{v}_H \, ,
\]
{\cred to be complemented with null initial conditions as in \eqref{3-2bis}.} 
We set $v = w_t(t)$ in the first equation and $v = u(t)$ in the second one, 
integrate over $(0, t)$ and add the two equations; it is straightforward to obtain
\[
\int_0^t \nh{w_t(s)} ds + \frac{\alpha}{2}\nh{\nabla w_t(t)} + \frac{1}{2}\nh{u(t)} + \int_0^t \nh{\nabla u(s)} ds \leq c\int_0^t \nh{u(s)} ds \, .
\]
According to the Gronwall lemma {\cred and owing to $w(0)=0$, it turns out that $w = u = 0$} a.e. in $Q$ and, by comparison in the second equation, $\xi = 0$ a.e. in $Q$.

\newcommand{\ww}{\widehat{w}_\beta}
\newcommand{\uu}{\widehat{u}_\beta}
\newcommand{\xx}{\widehat{\xi}_\beta}
\newcommand{\wz}{\widehat{w}_{0, \beta}}
\newcommand{\uz}{\widehat{u}_{0, \beta}}
\newcommand{\vz}{\widehat{v}_{0, \beta}}
\newcommand{\ff}{\widehat{f}_\beta}

\bigskip
\textbf{Error equations.} Because of the uniqueness, the whole family $\left\{(w_\beta, \, u_\beta, \, \xi_\beta)\right\}_{\beta > 0}$ {\cred converges}, as $\beta\searrow 0$, to the solution $(w, \, u, \, \xi)$ of Problem~$\Pa$. So, it makes sense to study the speed of this convergence. In order to perform that, we set $\ww = w_\beta - w$, $\uu = u_\beta - u$, $\xx = \xi_\beta - \xi$ and consider the problem obtained for these variables, by subtracting side by side the equations of Problems~$\Pab$ 
and~$\Pa$. For all $v\in V$ and a.a. {\cred $t\in (0,T)$, the equalities}
\begin{equation}
\begin{split}
\dual{\partial_t^2\ww(t)}{v} + \alpha\scal{\nabla \partial_t\ww(t)}{\nabla v}_H + \beta\scal{\nabla w_\beta(t)}{\nabla v}_H 
+ \dual{\partial_t \uu(t)}{v} \\ = \dual{\ff(t)}{v}
\end{split}
\label{eq A err}
\end{equation}
\begin{equation}
\begin{split}
\dual{\partial_t\uu(t)}{v} + \scal{\nabla \uu(t)}{\nabla v}_H + \scal{\xx(t)}{v}_H + \scal{g(u_\beta)(t) - g(u)(t)}{v}_H \\
= \scal{\partial_t\ww(t)}{v}_H 
\end{split}
\label{eq B err}
\end{equation}
are satisfied, as well as the {\cred initial conditions}
\[
\ww(0) = \wz \, , \qquad \partial_t\ww(0) = \vz \, , \qquad \uu(0) = \uz \, ,
\]
where {\cred $ \ff =  f_\beta - f = \ff^{(1)} + \ff^{(2)},$
\[
\begin{split}
\ff^{(1)}=  f_\beta^{(1)} - f^{(1)} \ \longrightarrow \ 0 \  \hbox{ in } \,
L^2 (0,T; V') \\ 
\ff^{(2)}=  f_\beta^{(2)} - f^{(2)}  \ \longrightarrow \ 0 \  \hbox{ in } \,  
L^1 (0,T; H)
\end{split}
\]
(cf.~\eqref{f rate a}),} $\wz := w_{0, \beta} - w_0$, $\vz := v_{0, \beta} - v_0$, and $\uz := u_{0, \beta} - u_0$.

\bigskip
\textbf{First estimate for the convergence error.} Now, we want to show Theorem~\ref{th: first estimate error}, so we assume all the needed hypotheses. Choose $v = \uu(t)$ in the equation \eqref{eq B err} and integrate over $(0, t)$; by the monotonicity of $\gamma$ and the Lipschitz-continuity of $g$, we easily derive
\begin{equation}
\frac{1}{2}\nh{\uu(t)} + \int_0^t \nh{\nabla \uu(s)} ds \leq  {\cred \frac{1}{2}\nh{\uz} + c \int_0^t \nh{\uu(s)} ds + \int_{Q_t} \uu\,\partial_t\ww \, }.
\label{test err 1, 1}
\end{equation}
We integrate with respect to time the equation \eqref{eq A err}:{\cred
\[
\begin{split}
\scal{\partial_t \ww(t)}{v}_H + \alpha \scal{\nabla \ww (t)}{\nabla v}_H + \beta \scal{1*\nabla w_\beta(t)}{\nabla v}_H + \scal{\uu(t)}{v}_H \\
= \langle 1*\ff(t), v \rangle + \scal{\vz + \uz}{v}_H + \alpha\scal{\nabla\wz}{\nabla v}_H \, .
\end{split}
\]
We set $v = \partial_t \ww$ and integrate over $(0, t)$; keeping only the first two terms in the left-hand side, we obtain
\begin{equation}
\begin{split}
\int_0^t \nh{\partial_t\ww(s)} ds + \frac{\alpha}{2}\nh{\nabla\ww(t)} \leq 
\frac{\alpha}{2}\nh{\nabla\wz} \\
- \beta \scal{1*\nabla w_\beta(t)}{\nabla \ww (t)}_H 
+ \beta\int_0^t\!\! \scal{\nabla w_\beta(s)}{\nabla \ww(s)}_H ds
-\int_{Q_t} \uu\,\partial_t\ww  \\
+ \int_0^t\!\! \left\langle 1*\ff^{(1)}(s) + \vz,
\partial_t \ww (s)\right\rangle  ds 
+  \int_{Q_t}\!\!\left(1* \ff^{(2)}+ \uz\right)\partial_t \ww + \cblu{\alpha\int_{Q_t} \nabla\wz \nabla\partial_t\ww} \, . 
\end{split}
\label{test err 1, 2}
\end{equation}
Due to the Young and H\"older inequalities and the boundedness of $\{w_\beta\}$ 
in $L^2(0,T;V)$, we have that 
\begin{equation}
\begin{split}
- \beta \scal{1*\nabla w_\beta(t)}{\nabla \ww (t)}_H 
\leq \frac{c}{\alpha} \beta^2 \int_0^t\!\! \nh{\nabla w_\beta (s)} ds 
+ {\cblu \frac{\alpha}{12}}\nh{\nabla\ww(t)} \\
\leq c\beta^2 + {\cblu \frac{\alpha}{12}}\nh{\nabla\ww(t)}
\end{split}
\label{nuova1}
\end{equation}
and
\begin{equation}
\beta\int_0^t\!\! \scal{\nabla w_\beta(s)}{\nabla \ww(s)}_H ds
\leq c\beta^2 +  \alpha \int_0^t\!\! \nh{\nabla\ww(s)}ds \, ,
\label{nuova2}
\end{equation}
{\cblu
\begin{equation}
\alpha \int_{Q_t}\nabla\wz \nabla \partial_t\ww \leq \frac{\alpha}{12}\nh{\nabla\ww(t)} + c\alpha\nh{\nabla\wz} \, .
\label{nuova5}
\end{equation} }
On the other hand, arguing as in the estimate of the term $T_4(t) $ 
of \eqref{test 1, 1} we deduce that
\begin{equation}
\begin{split}
\int_0^t\!\! \left\langle 1*\ff^{(1)}(s) + \vz,
\partial_t \ww (s)\right\rangle  ds \\
= \left\langle 1*\ff^{(1)}(t) + \vz, \ww (t )\right\rangle
- \int_0^t\!\! \left\langle \ff^{(1)}(s), \ww (s)\right\rangle  ds \\
\leq c \left( \int_0^t\norm{\ff^{(1)}(s)}_{V'}^2 ds + \Vert\vz \Vert_{V'}^2 
+ \Vert\wz \Vert_{H}^2  \right) + \frac14 \int_0^t\!\!  \nh{\partial_t\ww(s)} ds \\
+  {\cblu\frac{\alpha}{12}}\nh{\nabla\ww(t)} +  c  \int_0^t\!\! \left( \int_0^s \nh{\partial_t\ww(\tau)} d\tau \right) ds + 
c \, \alpha \int_0^t\!\! \nh{\nabla\ww(s)}ds 
\end{split}
\label{nuova3}
\end{equation}
Finally, we observe that 
\begin{equation}
\int_{Q_t}\!\!\left(1* \ff^{(2)}+ \uz\right)\partial_t \ww
\leq c \left( \norm{\ff^{(2)} }^2_{L^1(0,T;H)} + \Vert\uz \Vert_{H}^2  \right) + \frac14 \int_0^t\!\!  \nh{\partial_t\ww(s)} ds \\
\label{nuova4}
\end{equation}
Now we add \eqref{test err 1, 1} and \eqref{test err 1, 2}; collecting also
all the estimates in \eqref{nuova1}--\eqref{nuova4}, we find out that 
\[
\begin{split}
\frac{1}{2}\nh{\uu(t)} + \int_0^t \nh{\nabla \uu(s)} ds + \frac{1}{2}\int_0^t \nh{\partial_t\ww(s)} ds + \frac{\alpha}{4}\nh{\nabla\ww(t)} \\
\leq c\beta^2 + c \left( \norm{ \ff^{(1)} }_{L^2(0,T;V')}^2 
+ \norm{\ff^{(2)} }^2_{L^1(0,T;H)} + \Vert\uz \Vert_{H}^2 + \Vert\vz \Vert_{V'}^2 
  + \Vert\wz \Vert_{V}^2  \right) \\ 
+ c\int_0^t \nh{\uu(s)} ds 
+ c \int_0^t\!\! \left( \int_0^s \nh{\partial_t\ww(\tau)} d\tau \right) ds + 
c  \, \alpha \int_0^t\!\! \nh{\nabla\ww(s)}ds  .
\end{split}
\]
At this point, it suffices to recall \eqref{f rate a}--\eqref{data rate a}
and apply the Gronwall lemma to obtain the thesis of Theorem \ref{th: first estimate error}.
}%

\bigskip
\textbf{Second estimate for the convergence error.} 
{\cred Our aim is to prove Theorem \ref{th: second estimate error}, whose hypotheses are assumed to be satisfied. Thus, we can apply Theorems~\ref{th: regularity} and~\ref{th: L^infty estimate} to get a bound
\begin{equation}
\norm{u_\beta}_{L^\infty(Q)} + \norm{u}_{L^\infty(Q)} + \norm{\xi_\beta}_{L^\infty(Q)} + \norm{\xi}_{L^\infty(Q)}  \leq c_\alpha
\label{compactness}
\end{equation}
with $c_\alpha $ which is independent of $\beta$.
Now, if $\gamma$ is a maximal monotone graph which reduces 
to a single-valued  function in its domain, then $D(\gamma)$ is an open interval 
$(a, \, b)$ and, if $b < +\infty$, then $\gamma(r) \nearrow +\infty$ as 
$r \nearrow b$; similarly, if $a > -\infty$ then $\gamma(r) \searrow -\infty$ as 
$r\searrow a$. In any case, the condition \eqref{compactness} implies the 
existence of some compact interval $K\subseteq D(\gamma)$ such that 
$u_\beta(\overline Q)\subseteq K$ for all $\beta>0$, $u(\overline Q)\subseteq K$. 
Since $\gamma$ is assumed to be locally Lipschitz-continuous (cf.~\eqref{gamma 
lipshitz}), thanks to \eqref{stimaerr1} we immediately deduce that
\[
\norm{\xi_\beta - \xi}_{L^\infty(0,T; H)} \leq c 
\norm{u_\beta - u}_{L^\infty(0,T; H)}\leq c \beta \, .
\]
Moreover, by suitably modifying $g$ we can set $\xx \equiv 0$ in 
equation~\eqref{eq B err}, without loss of generality.

We start by taking $v = \partial_t\ww$ in \eqref{eq A err}, $v = \partial_t\uu$ in 
\eqref{eq B err}, integrating both equations over $(0, t)$ and adding side by 
side. Thanks to the Lipschitz-continuity of $g$ and the Young and H\"older inequalities, it is straightforward to obtain 
\[
\begin{split}
\frac{1}{2}\nh{\partial_t\ww(t)} + \alpha \int_0^t\!\!\nh{\nabla\partial_t\ww(s)}ds + 
\int_0^t\!\!\nh{\partial_t\uu(s)}ds +\frac{1}{2}\nh{\nabla\uu(t)} \\
\leq \frac{\beta^2}{\alpha} \norm{w_\beta }_{L^2(0,T;V)}^2 + 
\frac{\alpha}{4}\int_0^t\nh{\nabla \partial_t\ww(s)}ds 
+ c\int_0^t\nh{\uu(s)}ds 
+ \frac12\int_0^t\nh{\partial_t\uu(s)}ds \\
+  \frac1\alpha  \norm{\ff^{(1)}}^2_{L^2(0,T;V')} + 
\frac{\alpha}4 \norm{\partial_t\ww}^2_{L^2(0,T;H)} +
\frac{\alpha}{4}\int_0^t\!\!\nh{\nabla\partial_t\ww(s)}ds\\
+ \int_0^t\!\!\norm{\ff^{(2)}(s)}_H \norm{\partial_t\ww(s)}_H ds
+ \frac{1}{2}\nh{\vz} + 
\frac{1}{2}\nh{\nabla\uz} \, .
\end{split}
\]
Taking into account conditions \eqref{f rate a}, \eqref{data rate strong}
and the previous estimate \eqref{stimaerr1}, we easily have 
\[
\begin{split}
\frac{1}{2}\nh{\partial_t\ww(t)} + \frac{\alpha}2
 \int_0^t\nh{\nabla\partial_t\ww(s)}ds + \frac{1}{2} \int_0^t\nh{\partial_t\uu(s)}ds 
+\frac{1}{2}\nh{\nabla\uu(t)} \\
\leq c\, \beta^2 + \int_0^t\norm{\ff^{(2)} (s)}_H\norm{\partial_t\ww(s)}_H ds
\end{split}
\]
whence, by \eqref{f rate a} and a generalised Gronwall lemma {\cred (cf., e.g.,  
\cite[Lemme~A5, p.~157]{Brezis})}, we infer that
\begin{equation}
\norm{\ww}_{ \vett{W^{1, \, \infty}}{H} \cap \vett{H^1}{V}} +
\norm{\uu}_{ \vett{H^1}{H} \cap \vett{L^\infty}{V}} \leq c\, \beta
\label{estimate err 2}
\end{equation}
where the constant $c$ obviously depends on $\alpha$. 

Next, observe that the assumptions on the data are strong enough to guarantee that \eqref{eq A err} and \eqref{eq B err} can be reformulated as 
\begin{equation}
\partial_t^2 \ww - \alpha \Delta \partial_t\ww = \beta\Delta w_\beta - \partial_t \uu + \ff   \quad \hbox{ a.e. in } \, Q
\label{equazA}
\end{equation}
\begin{equation}
\partial_t\uu - \Delta \uu + g(u_\beta ) - g(u) = \partial_t \ww 
\quad \hbox{ a.e. in } \, Q
\label{equazB}
\end{equation}
along with the homogeneous Neumann boundary conditions for both $\ww$ and $\uu$. 

In view of \eqref{estimate err 2}, by a comparison of terms in \eqref{equazB} 
it is standard to deduce that $\norm{\Delta \uu}_{L^2(0,T;H)} \leq c\, \beta$ and consequently, 
owing to elliptic regularity estimates, we obtain
 \begin{equation}
\norm{\uu}_{\vett{L^2}{W}} \leq c_\alpha\beta \, .
\label{estimate err 2, 3}
\end{equation}
At this point, let us emphasize that for the proof of \eqref{estimate err 2} and 
\eqref{estimate err 2, 3} we have just used the control \eqref{f rate a}
on the difference $\ff$. 

We now pay attention to the equation \eqref{equazA} and 
multiply both sides by $-\Delta\partial_t\ww $, which belongs to $L^2(Q)$ 
(cf.~\eqref{w strong}), and integrate, also by parts, over $Q_t$. By means of the H\"older and Young inequalities, we infer that
\[
\begin{split}
\frac{1}{2}\nh{\nabla\partial_t\ww(t)} + \alpha \int_0^t\nh{\Delta\partial_t\ww(s)}ds \leq  \frac{1}{2}\nh{\nabla\vz} +
\frac{\beta^2}{\alpha} \int_0^t\nh{\Delta w_\beta(s)}ds\\
+ \frac{2}{\alpha}\int_0^t\nh{\partial_t\uu(s)}ds 
+ \frac{2}{\alpha} \norm{\ff^{(1)}}^2_{L^2(0,T;H)} 
+ \frac{\alpha}{2}\int_0^t\nh{\Delta\partial_t\ww(s)}ds \\
+ \int_0^t \norm{\nabla \ff^{(2)} (s) }_H
\norm{\nabla\partial_t\ww(s)}_H ds .
\end{split}
\]
Hence, recalling the uniform boundedness of $\{w_\beta\}$ in $L^2(0,T; W)$, we use  \eqref{data rate strong}, \eqref{estimate err 2}, \eqref{f rate b} and apply  
the generalised Gronwall lemma as before to obtain 
\[
\nh{\nabla\partial_t\ww(t)} + \int_0^t\nh{\Delta\partial_t\ww(s)}ds \leq c\, \beta^2.
\]
Now, by virtue of \eqref{fond. t. calculus} and \eqref{data rate strong} we also infer 
\[
\norm{\Delta \ww (t)}_H  \leq c\, \beta \quad \hbox{ for all } \, t\in [0,T] . 
\]
Then, standard elliptic regularity properties and the previous estimates 
\eqref{estimate err 2} and \eqref{estimate err 2, 3} lead us to \eqref{stimaerr2},
thus completing the proof of Theorem~\ref{th: second estimate error}. 



}%
\end{document}